\documentclass[10pt]{article}

\usepackage{amsmath}
\usepackage{mathtools}
\usepackage{stackrel}
\usepackage{scalerel}
\usepackage{leftidx}
\usepackage{xfrac}
\usepackage{stmaryrd}

\usepackage[charter]{mathdesign}
\usepackage{euscript}
\let\mathbb\undefined
\usepackage{bbold}
\DeclareSymbolFont{usualmathcal}{OMS}{cmsy}{m}{n}
\DeclareSymbolFontAlphabet{\mathcal}{usualmathcal}

\usepackage{authblk}

\usepackage[shortlabels]{enumitem}

\usepackage{hyperref}
\usepackage{cleveref}
\hypersetup{colorlinks=true, linkcolor=red!50!black, citecolor=green!50!black, urlcolor=blue!80!black}

\usepackage{xcolor}
\usepackage{tcolorbox}
\usepackage{graphicx}
\usepackage{subfigure}

\usepackage[all,2cell,arrow,matrix]{xy} \UseAllTwocells \SilentMatrices
\usepackage{tikz-cd}
\usepackage{tikz}
\usetikzlibrary{arrows,arrows.meta,%
decorations.markings,decorations.pathreplacing,%
decorations.pathmorphing,calc}

\usepackage{verbatim}
\usepackage{comment}

\usepackage[nottoc]{tocbibind} 
\usepackage{appendix}

\usepackage{geometry}
\tikzset{->-/.style={decoration={markings,mark=at position #1 with {\arrow{Stealth}}},postaction={decorate}},->-/.default=0.55}

\usepackage[amsmath,amsthm,thmmarks]{ntheorem}
\theoremstyle{definition}

\newtheorem{thm}{Theorem}[section]
\newtheorem{prop}[thm]{Proposition}

\newtheorem{cor}[thm]{Corollary}
\newtheorem{lem}[thm]{Lemma}

{
\theoremsymbol{\mbox{${\blacksquare}$}}
\newtheorem{defn}[thm]{Definition}
}
{
\theoremsymbol{\mbox{$\heartsuit$}}
\newtheorem{expl}[thm]{Example}
}
{
\theoremsymbol{\mbox{$\diamondsuit$}}
\newtheorem{rem}[thm]{Remark}
}
\qedsymbol{\mbox{$\square$}}

\numberwithin{equation}{section}
\numberwithin{thm}{section}

\newcommand\be            {\begin{equation}}
\newcommand\ee            {\end{equation}}
\newcommand\bea           {\begin{eqnarray}}
\newcommand\eea         {\end{eqnarray}}
\newcommand\bnu          {\begin{enumerate}}
\newcommand\enu          {\end{enumerate}}
\newcommand\bit          {\begin{itemize}}
\newcommand\eit          {\end{itemize}}

\newcommand{\pf}{\begin{proof}}
\newcommand{\epf}{\end{proof}}

\makeatletter
\providecommand{\leftsquigarrow}{%
  \mathrel{\mathpalette\reflect@squig\relax}%
}
\newcommand{\reflect@squig}[2]{%
  \reflectbox{$\m@th#1\rightsquigarrow$}%
}
\makeatother

\newcommand\Nb			{\mathbb{N}}

\newcommand\Zb			{\mathbb{Z}}

\newcommand\bk			{\mathbb{k}}

\newcommand\CA			{\EuScript{A}}
\newcommand\CB			{\EuScript{B}}
\newcommand\CC			{\EuScript{C}}
\newcommand\CD			{\EuScript{D}}
\newcommand\CE			{\EuScript{E}}

\newcommand\CG			{\EuScript{G}}

\newcommand\CM			{\EuScript{M}}
\newcommand\CN			{\EuScript{N}}

\newcommand\CV			{\EuScript{V}}
\newcommand\CW			{\EuScript{W}}

\newcommand{\FZ}			{\text{\usefont{U}{euf}{m}{n}Z}}

\newcommand\SC			{\mathsf{C}}
\newcommand\SD			{\mathsf{D}}

\DeclareMathOperator{\Hom}{Hom}
\DeclareMathOperator{\End}{End}
\DeclareMathOperator{\Aut}{Aut}
\DeclareMathOperator{\Ext}{Ext}

\DeclareMathOperator{\id}{id}

\DeclareMathOperator{\ob}{ob}

\DeclareMathOperator{\fun}{Fun}

\DeclareMathOperator{\Idem}{Idem}

\DeclareMathOperator{\LMod}{LMod}
\DeclareMathOperator{\RMod}{RMod}

\DeclareMathOperator{\Stab}{Stab}

\newcommand{\op}			{\mathrm{op}}
\newcommand{\rev}			{\mathrm{rev}}

\newcommand{\one}			{\mathbb{1}}

\newcommand\set			{\mathrm{Set}}

\newcommand\vect			{\mathrm{Vec}}
\newcommand\svect			{\mathrm{sVec}}

\newcommand\rep			{\mathrm{Rep}}

\newcommand\Irr			{\mathrm{Irr}}

\newcommand{\bscale}	{0.7}
\makeatletter
\newcommand{\ec}[2][]	{{\@ec{#1 |}{#2}}}
\newcommand{\bc}[2][]	{{\@ec{#1}{#2}}}
\newcommand{\@ec}[2]	{\mathchoice
  {\displaystyle \raise.9ex\hbox{$\scaleobj{\bscale}{#1}$} {#2}}%
  {\textstyle \raise.9ex\hbox{$\scaleobj{\bscale}{#1}$} {#2}}%
  {\scriptstyle \raise.55ex\hbox{$\scriptstyle \scaleobj{\bscale}{#1}$} {#2}}%
  {\scriptscriptstyle \raise.38ex\hbox{$\scriptscriptstyle \scaleobj{\bscale}{#1}$} {#2}}%
}
\makeatother

\DeclareMathOperator{\Kar}{Kar}

\newcommand{\EE}{\mathrm{E}}

\begin{document}

\title{\huge Tannaka-Krein duality for finite 2-groups}
\author[a]{Mo Huang \thanks{Email: \href{mailto:jasmine.huang.527@gmail.com}{\tt jasmine.huang.527@gmail.com}}}
\author[c,d]{Zhi-Hao Zhang \thanks{Email: \href{mailto:zhangzhihao@bimsa.cn}{\tt zhangzhihao@bimsa.cn}}}
\affil[a]{School of Mathematical Sciences, East China Normal University, Shanghai, 200241, China}
\affil[b]{Beijing Institute of Mathematical Sciences and Applications, Beijing, 101408, China}
\affil[c]{Wu Wen-Tsun Key Laboratory of Mathematics of Chinese Academy of Sciences, \authorcr
School of Mathematical Sciences, University of Science and Technology of China, Hefei, 230026, China}\date{\vspace{-5ex}}

\maketitle

\begin{abstract}
Let $\CG$ be a finite 2-group. We show that the 2-category $2\rep(\CG)$ of finite semisimple 2-representations is a symmetric fusion 2-category. We also relate the auto-equivalence 2-group of the symmetric monoidal forgetful 2-functor $\omega \colon 2\rep(\CG) \to 2\vect$ to the auto-equivalence 2-group of the regular algebra and show that they are equivalent to $\CG$. This result categorifies the usual Tannaka-Krein duality for finite groups.
\end{abstract}

\tableofcontents

\section{Introduction}

For a finite group $G$, the well-known Tannaka-Krein duality \cite{Tan39,Kre49} says that the automorphism group $\Aut^{\EE_3}(\omega)$ of symmetric monoidal natural automorphisms of the forgetful functor $\omega \colon \rep(G) \to \vect$ is canonically isomorphic to $G$ (see also \cite[Section 1]{JS91a} for a quick review). In other words, the group $G$ can be reconstructed from the category $\rep(G)$ of finite-dimensional $G$-representations.

It is natural to categorify the Tannaka-Krein duality to (finite) 2-groups. A 2-group is a categorification of a group. There are three equivalent models for 2-group: weak 2-group, strict 2-group and crossed module (of groups). A weak 2-group is a monoidal category whose morphisms and objects are invertible. A strict 2-group is a weak 2-group whose objects are invertible in a strict sense. A crossed module consists of two groups $G,H$, a group homomorphism $G \to H$ and an $H$-action on $G$, satisfying certain axioms. It is tautological that the notion of a strict 2-group and a crossed module are equivalent to the notion of a category internal to the category of groups (see for example \cite[Section 7]{BL04}). Moreover, the notion of a weak 2-group is also equivalent to that of a strict 2-group \cite[Proposition 45]{BL04}. In this work, we use the weak 2-group model because we usually consider skeletal 2-groups, which are typically not strict.

To establish the Tannaka-Krein duality for finite 2-groups, the first step is to define the 2-category of 2-representations of a 2-group. The main issue is to find the 2-category of `2-vector spaces' in which we represent 2-groups. There are two commonly used definitions of a 2-vector space: the Kapranov-Voevodsky 2-vector space \cite{KV94} and the Baez-Crans 2-vector space \cite{BC04}. A Kapranov-Voevodsky 2-vector space is a finite direct sum of the category $\vect$ of finite-dimensional vector spaces. When the ground field $\bk$ is algebraically closed, it is also equivalent to a finite semisimple category. A Baez-Crans 2-vector space is a category internal to $\vect$, or equivalently, a 2-term chain complex of vector spaces. There were several works towards the 2-representation theory of 2-groups in either the 2-category of Kapranov-Voevodsky 2-vector spaces \cite{BM06,Elg07,Elg11} or the 2-category of Baez-Crans 2-vector spaces \cite{FB03,HE16}. More generally, there were also works on `infinite-dimensional' 2-representations \cite{CY05,BBFW12} based on the notion of a measurable category \cite{Yet05}.

In this work, we are mainly interested in `finite-dimensional 2-representations'. Recall that the dualizable objects in the monoidal category of vector spaces are exactly finite-dimensional ones. This motivates us to use a fully dualizable 2-category to model the 2-category of `finite-dimensional 2-vector spaces'. A candidate is the 2-category of separable algebras, finitely generated projective bimodules and bimodule homomorphisms. When the ground field $\bk$ is perfect, this 2-category is also equivalent to the 2-category $2\vect$ of finite semisimple categories, linear functors and natural transformations. When the ground field $\bk$ is algebraically closed, $2\vect$ is also equivalent to the 2-category of Kapranov-Voevodsky 2-vector spaces. We call a 2-representation of a 2-group $\CG$ in $2\vect$ a \emph{finite semisimple 2-representation}, and the 2-category of finite semisimple 2-representations is denoted by $2\rep(\CG)$.

A basic result in the representation theory of groups is Maschke's theorem \cite{Mas98,Mas99}, which states that the category $\rep(G)$ of finite-dimensional representations of a finite group $G$ is finite semisimple if the characteristic of the ground field does not divide the order of $G$. It is natural to ask whether $2\rep(\CG)$ for a finite 2-group $\CG$ is finite semisimple, in the sense of Douglas and Reutter \cite{DR18} when the ground field is algebraically closed of characteristic zero\footnote{Over a general ground field, the notion of a finite semisimple 2-category should be replaced by that of a \emph{compact semisimple 2-category} \cite{Dec23}.}. The proof has been sketched in \cite[Example 1.4.19]{DR18}. The main idea is to categorify the fact that the representations of a group $G$ are equivalent to the modules over the group algebra $\bk[G]$. In the same spirit, we define a `2-group 2-algebra' $\vect_\CG$ for a finite 2-group $\CG$ and show that $2\rep(\CG)$ is equivalent to the 2-category of finite semisimple left modules over $\vect_\CG$. Then the semisimplicity of $2\rep(\CG)$ follows from the semisimplicity of $\vect_\CG$. Also, the structure of $\vect_\CG$ helps us to analyze the structure of $2\rep(\CG)$, including the simple objects, hom categories, fusion rules, etc (see Theorem \ref{thm_2Rep_G_structure}). 

The dual of $\vect_\CG$ is naturally an algebra in $2\rep(\CG)$ and called the regular representation \cite{Elg11} or the regular algebra of $\CG$. It is the category $\fun(\CG,\vect)$ of functors from $\CG$ to $\vect$, categorifying the function algebra $\fun(G)$ for finite groups $G$. The regular algebra also plays an important role in the proof of the Tannaka-Krein duality for finite 2-groups. If we forget the $\CG$-action on $\fun(\CG,\vect)$, we can also study the finite semisimple modules over this multi-fusion category. It turns out that these modules are equipped with a 2-group $\CG$ grading, which generalize the notion of a group-graded category. 

Both the 2-categories of modules over $\vect_\CG$ and $\fun(\CG,\vect)$ have natural monoidal structures. These monoidal structures are naturally induced from the Hopf monoidal structure \cite{CF94,Neu97} on $\vect_\CG$ and $\fun(\CG,\vect)$, respectively. By the Tannaka-Krein reconstruction for fusion 2-categories \cite{Gre23}, these Hopf monoidal structures can also be reconstructed from the forgetful 2-functors $2\rep(\CG) \to 2\vect$ and $2\vect_\CG \to 2\vect$, respectively. We explicitly compute the Hopf monoidal structures on $\vect_\CG$ and $\fun(\CG,\vect)$.

Our result on the Tannaka-Krein duality for finite 2-groups is the first step towards the classification of symmetric fusion 2-categories. For symmetric fusion 1-categories, the classification is given by Deligne's theorem \cite{Del07,Del02}:
\bit
\item Every symmetric fusion category $\CE$ admits a symmetric fiber functor $\omega \colon \CE \to \svect$. It is unique up to natural isomorphism.
\item Let $G$ be the group $\Aut^{\EE_3}(\omega)$ of symmetric natural automorphisms of $\omega$. It is known that $\Aut^{\EE_3}(\id_{\svect}) \simeq \Zb_2$, where the generator is the parity morphisms $\Pi = \{\Pi_V \colon V \to V\}_{V \in \svect}$ defined by $\Pi_V = \id_{V_0} \oplus (-\id_{V_1})$. Let $z \in G$ be the whiskering of $\Pi$ and $\omega$. It is easy to see that $z$ is central and $z^2 = e$. We say that $(G,z)$ is a super group. Then $\CE$ is equivalent to the symmetric monoidal category $\mathrm{sRep}(G,z)$ of finite-dimensional super representations of $(G,z)$.
\item In particular, if the symmetric fiber functor $\omega$ factors through $\vect \hookrightarrow \svect$ (called the `bosonic case'), then $z \in G$ is trivial and $\Aut^{\EE_3}(\omega \colon \CE \to \svect)$ is the same as the group $\Aut^{\EE_3}(\CE \to \vect)$. In this case, $\CE$ is equivalent to $\rep(G)$.
\eit
It is natural to generalize this theorem to the classification of symmetric fusion 2-categories. Recently, D\'{e}coppet and Yu showed in \cite{DY23} that every symmetric fusion 2-category admits a symmetric fiber 2-functor to $2\svect$. Then the next step should be to establish the super 2-representation theory of super 2-groups. Our result in this paper is a justification of the bosonic case (that the super structure is trivial): the autoequivalence 2-group of the forgetful 2-functor $2\rep(\CG) \to 2\vect$ is equivalent to $\CG$.

\smallskip
The outline of this paper is as follows. In Section \ref{sec:preliminary} we review some basic notions and facts that are useful in the study of 2-groups and 2-representations. In Section \ref{sec:2-group}, first we review the basic 2-group theory, and then study the 2-category $2\rep(\CG)$ of finite semisimple 2-representations of a finite 2-group $\CG$. To show that it is a finite semisimple 2-category, we construct the 2-group 2-algebra $\vect_\CG$ and show that $2\rep(\CG)$ is equivalent to the 2-category of finite semisimple modules over $\vect_\CG$. We also study the dual $\fun(\CG,\vect)$ of $\vect_\CG$ and its modules. In particular, we explicitly compute the Hopf monoidal structure on $\vect_\CG$ and $\fun(\CG,\vect)$. Finally, in Section \ref{sec:duality} we establish the Tannaka-Krein duality for finite 2-groups.

\medskip
\noindent \textbf{Notations and Conventions.} For a group $G$, we usually denote its unit by $e \in G$. Throughout this paper, $\bk$ is an algebraically closed field of characteristic zero. We use $\vect$ to denote the $\bk$-linear category of finite dimensional $\bk$-vector spaces, and use $\rep(G)$ to denote the category of finite-dimensional representations of $G$. For a finite abelian group $A$, its dual group, i.e., the group of characters, is denoted by $\hat A$. The identity morphism of an object $x$ in a category $\CC$ is denoted by $1_x$ or simple $1$. The opposite category of a category $\CC$ obtained by reversing the morphisms is denoted by $\CC^\op$. For a monoidal category $\CD$, the same underlying category $\CD$ equipped with the reversed tensor product is denoted by $\CD^\rev$. By a 2-category, a 2-functor and a 2-natural transformation, we mean the notions that are usually called a weak 2-category (or a bicategory), a pseudofunctor and a strong transformation, respectively. Given two categories (2-categories, resp.) $X,Y$, we use $\fun(X,Y)$ to denote the category (2-category, resp.) of functors (2-functors, resp.) from $X$ to $Y$; if $X$ and $Y$ are $\bk$-linear, $\fun_\bk(X,Y)$ denotes the $\bk$-linear category ($\bk$-linear 2-category, resp.) of $\bk$-linear functors ($\bk$-linear 2-functors, resp.) from $X$ to $Y$.

\bigskip
\noindent \textbf{Acknowledgments.} We would like to thank Liang Kong and Hao Xu for helpful discussions. MH was supported by Research Grants Council of Hong Kong under GRF~17311322. ZHZ was supported by Wu Wen-Tsun Key Laboratory of Mathematics at USTC of Chinese Academy of Sciences and the start-up grant of BIMSA and the China Postdoctoral Science Foundation under Grant Number 2025M783072. MH and ZHZ was also supported by NSFC (Grant No.~11971219) and by Guangdong Provincial Key Laboratory (Grant No.~2019B121203002) and by Guangdong Basic and Applied Basic Research Foundation (Grant No.~2020B1515120100).

\section{Preliminaries} \label{sec:preliminary}

In this section we briefly review some basic notions and facts.

\subsection{Rigid monoidal categories}

Let $\CC$ be a monoidal category. For convenience, we omit the associator and unitor of $\CC$ in this subsection. By MacLane's strictness theorem \cite{Mac63,JS93}, we may also assume without loss of generality that $\CC$ is a strict monoidal category.

A \emph{left dual} of an object $x \in \CC$ is a triple $(x^L \in \CC,b_x \colon \one \to x \otimes x^L,d_x \colon x^L \otimes x \to \one)$ satisfying the following zig-zag equations:
\begin{gather*}
1_x = \bigl(x \simeq \one \otimes x \xrightarrow{b_x \otimes 1} x \otimes x^L \otimes x \xrightarrow{1 \otimes d_x} x \otimes \one \simeq x \bigr) , \\
1_{x^L} = \bigl(x^L \simeq x^L \otimes \one \xrightarrow{1 \otimes b_x} x^L \otimes x \otimes x^L \xrightarrow{d_x \otimes 1} \one \otimes x^L \simeq x^L \bigr) .
\end{gather*}
Dually, a \emph{right dual} of $x \in \CC$ is a triple $(x^R,b_x' \colon \one \to x^R \otimes x,d_x' \colon x \otimes x^R \to \one)$ satisfying similar equations. A left or right dual of $x$, if exists, is unique up to a unique isomorphism. A \emph{rigid monoidal category} is a monoidal category in which every object has a left and right dual.

Given a rigid monoidal category $\CC$, by fixing a left dual $(x^L,b_x,d_x)$ for every object $x \in \CC$, we can define a left dual functor $\delta^L \colon \CC^\op \to \CC$ which maps an object $x$ to its left dual $x^L$ and maps a morphism $f \colon x \to y$ to
\[
f^L \coloneqq \bigl( y^L \simeq y^L \otimes \one \xrightarrow{1 \otimes b_x} y^L \otimes x \otimes x^L \xrightarrow{1 \otimes f \otimes 1} y^L \otimes y \otimes x^L \xrightarrow{d_y \otimes 1} \one \otimes x^L \simeq x^L \bigr) .
\]
This functor can be promoted to a monoidal functor $\delta^L \colon \CC^\op \to \CC^\rev$. Different choices of left dual functors are canonically isomorphic.

\begin{lem} \label{lem:invertible_object}
Let $\CC$ be a monoidal category and $x \in \CC$. The following statements are equivalent:
\bnu[(1)]
\item There exist $y,z \in \CC$ such that $y \otimes x \simeq \one$ and $x \otimes z \simeq \one$.
\item There exists $y \in \CC$ such that $y \otimes x \simeq \one \simeq x \otimes y$.
\item There is a left dual $(x^L,b_x \colon \one \to x \otimes x^L,d_x \colon x^L \otimes x \to \one)$ of $x$ such that $b_x$ and $d_x$ are isomorphisms.
\item There is a right dual $(x^R,b'_x \colon \one \to x^R \otimes x,d'_x \colon x \otimes x^R \to \one)$ of $x$ such that $b'_x$ and $d'_x$ are isomorphisms.
\enu
\end{lem}

\pf
Clearly both (3) and (4) imply (1), and (1) is equivalent to (2). Now we show that (2) implies (3), and the proof of the implication (2) $\implies$ (4) is similar. Suppose $y \in \CC$ and $f \colon \one \to x \otimes y$, $g \colon y \otimes x \to \one$ are isomorphisms. Let
\begin{gather*}
u \coloneqq \bigl( x \simeq \one \otimes x \xrightarrow{f \otimes 1} x \otimes y \otimes x \xrightarrow{1 \otimes g} x \otimes \one \simeq x \bigr) , \\
v \coloneqq \bigl( y \simeq y \otimes \one \xrightarrow{1 \otimes f} y \otimes x \otimes y \xrightarrow{g \otimes 1} \one \otimes y \simeq y \bigr) .
\end{gather*}
In general, $u,v$ are isomorphisms but not the identity morphisms, and thus $(y,f,g)$ is not a left dual of $x$.

Since $f,g$ are isomorphisms, one can verify that both $1_y \otimes u$ and $v \otimes 1_x$ are equal to the composite morphisms
\[
y \otimes x \simeq y \otimes \one \otimes x \xrightarrow{1 \otimes f \otimes 1} y \otimes x \otimes y \otimes x \xrightarrow{g \otimes g} \one \xrightarrow{g^{-1}} y \otimes x .
\]
Then we define $b \coloneqq f$ and $d \coloneqq g \circ (v^{-1} \otimes 1_x) = g \circ (1_y \otimes u^{-1})$. It is not hard to see that $(y,b,d)$ is a left dual of $x$.
\epf

\begin{defn}
Let $\CC$ be a monoidal category. An object $x \in \CC$ is \emph{invertible} if it satisfies one (and hence all) of the conditions in Lemma \ref{lem:invertible_object}.
\end{defn}

It is clear that monoidal functors preserve invertible objects.

\subsection{Multi-fusion categories}

A \emph{multi-fusion category} is a $\bk$-linear rigid monoidal category that is finite semisimple, and a \emph{fusion category} is a multi-fusion category whose tensor unit is simple. Given two multi-fusion categories $\CC,\CD$, we say that they are \emph{Morita equivalent} if there exists a finite semisimple $\CC$-$\CD$-bimodule $\CM$ and a finite semisimple $\CD$-$\CC$-bimodule $\CN$ such that $\CM \boxtimes_\CD \CN \simeq \CC$ as $\CC$-$\CC$-bimodules and $\CN \boxtimes_\CC \CM \simeq \CD$ as $\CD$-$\CD$-bimodules. Such $\CM$ and $\CN$ are called the invertible bimodules. Equivalently, $\CC$ and $\CD$ are Morita equivalent if the 2-categories $\RMod_\CC(2\vect)$ and $\RMod_\CD(2\vect)$ are equivalent, where $\RMod_\CC(2\vect)$ is the 2-category of finite semisimple right $\CC$-modules, right $\CC$-module functors and right $\CC$-module natural transformations.

Let us recall the structure theorem of multi-fusion categories \cite[Section 4.3]{EGNO15} (see also \cite[Theorem 2.5.1]{KZ18}):

\begin{thm} \label{thm:structure_multi-fusion_category}
Let $\CC$ be a multi-fusion category and $\one = \bigoplus_{i \in \Lambda} \one_i$ be a decomposition of the tensor unit $\one \in \CC$ into simple objects. Denote $\CC_{i,j} \coloneqq \one_i \otimes \CC \otimes \one_j$ for $i,j \in \Lambda$.
\bnu[(1)]
\item The tensor product of $\CC$ induces functors $\CC_{i,j} \times \CC_{j,k} \to \CC_{i,k}$ for all $i,j,k \in \Lambda$. In particular, $\CC_{i,i}$ is a fusion category and $\CC_{i,j}$ is a finite semisimple $\CC_{i,i}$-$\CC_{j,j}$-bimodule for every $i,j \in \Lambda$.
\item The relation $\sim$ on $\Lambda$ defined by $i \sim j$ if $\CC_{i,j} \not\simeq 0$ is an equivalence relation.
\item For every equivalence class $x \subseteq \Lambda$, denote $\CC_{x,x} \coloneqq \bigoplus_{i,j \in x} \CC_{i,j}$. Then $\CC_{x,x}$ is an indecomposable multi-fusion category and $\CC \simeq \bigoplus_x \CC_{x,x}$ as multi-fusion categories.
\item For every equivalence class $x \subseteq \Lambda$ and $i,j,k \in x$, the functor $\CC_{i,j} \boxtimes_{\CC_{j,j}} \CC_{j,k} \simeq \CC_{i,k}$ induced by the tensor product of $\CC$ is an equivalence of $\CC_{i,i}$-$\CC_{k,k}$-bimodules. In particular, $\CC_{i,i}$ is Morita equivalent to $\CC_{j,j}$ and the invertible bimodules are $\CC_{i,j}$ and $\CC_{j,i}$. Moreover, $\CC_{i,i}$ is Morita equivalent to $\CC_{x,x}$.
\enu
\end{thm}

Theorem \ref{thm:structure_multi-fusion_category} can be viewed as a categorification of the classical Artin-Wedderburn theorem which states that every finite semisimple algebra is the direct sum of matrix algebras over division algebras.

\subsection{Karoubi completion} \label{sec_Kar}

Let $\CC$ be a category. A morphism $p \colon x \to x$ in $\CC$ is called an \emph{idempotent} if $p \circ p = p$. We say an idempotent $p \colon x \to x$ \emph{splits} if there exists an object $y \in \CC$ and morphisms $r \colon x \to y$ and $s \colon y \to x$ such that $s \circ r = p$ and $r \circ s = 1_y$. It is easy to see that $y$ is both the equalizer and coequalizer of $p$ and $1_x$, hence is unique up to a unique isomorphism.

A category is called \emph{idempotent complete} if every idempotent splits. An \emph{idempotent completion} of a category $\CC$ is an idempotent complete category $\Idem(\CC)$ equipped with a functor $\iota \colon \CC \to \Idem(\CC)$ satisfying the following universal property: for every idempotent complete category $\CD$, the functor
\[
- \circ \iota \colon \fun(\Idem(\CC),\CD) \to \fun(\CC,\CD)
\]
is an equivalence. A construction of $\mathrm{Idem}(\CC)$ is given as follows:
\bit
\item The objects are pairs $(x,p)$, where $x \in \CC$ and $p \colon x \to x$ is an idempotent.
\item A morphism $f \colon (x,p) \to (y,q)$ is a morphism $f \colon x \to y$ in $\CC$ such that $f \circ p = f = q \circ f$.
\eit
The natural embedding $\iota \colon \CC \to \Idem(\CC)$ sends $x \in \CC$ to $(x,1_x)$.

\begin{lem} \label{lem:idempotent_completion_monoidal}
Let $\CC$ be a monoidal category. Then $\Idem(\CC)$ has a natural monoidal structure and the natural embedding $\iota \colon \CC \to \Idem(\CC)$ is naturally a monoidal functor.
\end{lem}

\pf
The monoidal structure on $\Idem(\CC)$ can be induced by the universal property. Here we give an explicit construction.
\bit
\item The tensor product is given by $(x,p) \otimes (y,q) \coloneqq (x \otimes y,p \otimes q)$.
\item The tensor unit is $(\one,1_\one)$.
\item The associator $\alpha_{(x,p),(y,q),(z,r)} \colon ((x,p) \otimes (y,q)) \otimes (z,r) \to (x,p) \otimes ((y,q) \otimes (z,r))$ is the composite morphism
\[
(p \otimes (q \otimes r)) \circ \alpha_{x,y,z} \circ ((p \otimes q) \otimes r) .
\]
By the naturality of $\alpha$, this is also equal to $\alpha_{x,y,z} \circ ((p \otimes q) \otimes r)$ and $(p \otimes (q \otimes r)) \circ \alpha_{x,y,z}$.
\item The left unitor $\lambda_{(x,p)} \colon (\one,1_\one) \otimes (x,p) \to (x,p)$ is the composite morphism
\[
p \circ \lambda_x \circ (1_\one \otimes p) .
\]
By the naturality of $\lambda$, this is also equal to $\lambda_x \circ (1_\one \otimes p)$ and $p \circ \lambda_x$.
\item The right unitor $\rho_{(x,p)} \colon (x,p) \otimes (\one,1_\one) \to (x,p)$ is the composite morphism
\[
p \circ \rho_x \circ (p \otimes 1_\one) .
\]
By the naturality of $\rho$, this is also equal to $\rho_x \circ (p \otimes 1_\one)$ and $p \circ \rho_x$.
\eit
Then it is easy to verify that $\Idem(\CC)$ is a monoidal category and $\iota \colon \CC \to \Idem(\CC)$ is a monoidal functor with the trivial monoidal structure.
\epf

A $\bk$-linear category is called \emph{Karoubi complete} if it is additive (i.e., has all finite direct sums) and idempotent complete. A \emph{Karoubi completion} of a $\bk$-linear category $\CC$ is a Karoubi complete $\bk$-linear category $\Kar(\CC)$ equipped with a $\bk$-linear functor $\iota \colon \CC \to \Kar(\CC)$ satisfying the following universal property: for every Karoubi complete $\bk$-linear category $\CD$, the functor
\[
- \circ \iota \colon \fun_\bk(\Kar(\CC),\CD) \to \fun_\bk(\CC,\CD)
\]
is an equivalence. A construction of $\Kar(\CC)$ is given by first taking the additive completion (i.e., formally adding all finite direct sums) and then taking the idempotent completion.

\begin{expl} \label{expl:Karoubi_completion_group}
Let $R$ be a $\bk$-algebra. Its one-point delooping $\mathrm B R$ is a $\bk$-linear category. The Karoubi completion $\Kar(\mathrm B R)$ is equivalent to the category of finitely generated projective right $R$-modules. In particular, when $R = \bk[A]$ where $A$ is a finite abelian group, the Karoubi completion $\Kar(\mathrm B \bk[A])$ is equivalent to the category $\rep(A)$ of finite-dimensional representations of $A$.
\end{expl}

\subsection{Finite semisimple and multi-fusion 2-categories}

\begin{defn}[\cite{DR18}]
Let $\SC$ be a 2-category. A monad $(T \colon x \to x,\mu \colon T \circ T \Rightarrow T,\eta \colon 1_x \Rightarrow T)$ is \emph{separable} if there exists a $T$-$T$-bimodule map $\delta \colon T \Rightarrow T \circ T$ such that $\mu \cdot \delta = 1_T$. An adjunction $(f \colon x \to y,g \colon y \to x,\eta \colon 1_x \Rightarrow g \circ f,\varepsilon \colon f \circ g \Rightarrow 1_y)$ is \emph{separable} if there exists $\phi \colon 1_y \Rightarrow f \circ g$ such that $\varepsilon \cdot \phi = 1_{1_y}$. A separable monad $T$ \emph{separably splits} if there exists a separable adjunction $f \dashv g$ such that $T \simeq g \circ f$. A 2-category is \emph{idempotent complete} if every hom category is idempotent complete and every separable monad separably splits.
\end{defn}

The separable monads can be viewed as categorified idempotents.

\begin{rem}
Gaiotto and Johnson-Freyd further generalize the notion of idempotents and idempotent completeness to $n$-categories in \cite{GJF19}. A separable monad with a choice of the right inverse of the multiplication is a \emph{unital 2-condensation monad}, and a separable adjunction with a choice of the right inverse of the counit is a \emph{unital 2-condensation} as defined in \cite{GJF19}.
\end{rem}

\begin{defn}[\cite{DR18}] \label{defn:semisimple_2-cat}
A \emph{finite semisimple} 2-category is a $\bk$-linear 2-category $\SC$ satisfying the following conditions:
\bnu[(1)]
\item For every $x,y \in \SC$, the hom category $\Hom_\SC(x,y)$ is a finite semisimple category.
\item The direct sums of finitely many objects exist in $\SC$. In particular, $\SC$ has a zero object.
\item Every 1-morphism in $\SC$ has a left and right adjoint.
\item $\SC$ is idempotent complete.
\item There are finitely many equivalence classes of simple objects.
\enu
An object $x$ in a finite semisimple 2-category $\SC$ is called \emph{simple} if the identity 1-morphism $1_x \in \Hom_\SC(x,x)$ is a simple object.
\end{defn}

\begin{thm}[Categorical Schur's lemma, \cite{DR18}]
Let $\SC$ be a finite semisimple 2-category and $x,y,z \in \SC$ are simple objects. If $f \colon x \to y$ and $g \colon y \to z$ are non-zero 1-morphisms, then $g \circ f$ is also non-zero.
\end{thm}

In a finite semisimple 2-category $\SC$, two simple objects $x,y \in \SC$ are \emph{connected} if there exists a nonzero 1-morphism between them. By categorical Schur's lemma, being connected is an equivalence relation. An equivalence class is called a \emph{connected component} of $\SC$, and the set of connected components is denoted by $\pi_0(\SC)$.

\begin{thm}[\cite{DR18}] \label{thm:semisimple_2-cat_multi-fusion}
Let $\CC$ be a multi-fusion category. The 2-category $\LMod_\CC(2\vect)$ of finite semisimple left $\CC$-modules is a finite semisimple 2-category. Conversely, every finite semisimple 2-category is equivalent to $\LMod_\CA(2\vect)$ for some multi-fusion category $\CA$.
\end{thm}

\begin{rem}
If we work over an arbitrary field, the 2-category of finite semisimple (or even separable) right modules over a multi-fusion category may not be finite semisimple as defined in Definition \ref{defn:semisimple_2-cat} because there may be infinitely many simple objects up to equivalence. A better notion is the \emph{compact semisimple 2-category} defined by D\'{e}coppet \cite{Dec23}, with the finiteness condition being replaced by that the set of connected components is finite.
\end{rem}

\begin{defn}[\cite{DR18}]
A \emph{multi-fusion 2-category} is a finite semisimple monoidal 2-category such that every object has a left and right dual. A \emph{fusion 2-category} is a multi-fusion 2-category such that the tensor unit is simple.
\end{defn}

\subsection{Tannaka-Krein duality for finite groups} \label{sec:TK_dual_finite_group}

In this subsection we briefly review the proof of the Tannaka-Krein duality for finite groups. Let $G$ be a finite group and $\omega \colon \rep(G) \to \vect$ be the forgetful functor.

\begin{thm}
The group $\Aut^{\EE_3}(\omega)$ of symmetric monoidal auto-equivalences of $\omega$ is isomorphic to $G$.
\end{thm}

\pf
Given $g \in G$ there is a symmetric monoidal natural isomorphism $\Phi^g \colon \omega \Rightarrow \omega$ defined by $\Phi^g_{(V,\rho)} \coloneqq \rho(g) \colon \omega(V) \to \omega(V)$ for every $(V,\rho) \in \rep(G)$. The assignment $g \mapsto \Phi^g$ defines a group homomorphism $\Phi \colon G \to \Aut^{\EE_3}(\omega)$.

The construction of the inverse of $\Phi$ needs several observations.
\bit
\item Recall the vector space $\Hom(G,\bk)$ equipped with the $G$-action induced by right translation and the point-wise multiplication is a commutative algebra in $\rep(G)$, called the \emph{regular algebra} of $G$. We denote this algebra by $\fun^r(G)$. Similarly, the vector space $\Hom(G,\bk)$ equipped with the $G$-action induced by left translation and the point-wise multiplication is also a commutative algebra in $\rep(G)$, denoted by $\fun^l(G) \in \rep(G)$. We also call $\fun^l(G)$ the regular algebra. They are indeed isomorphic to each other.
\item The endomorphism algebra $\End_{\rep(G)}(\fun^r(G))$ is isomorphic to the group algebra $\bk[G]$, which acts on the vector space $\omega(\fun^r(G)) = \Hom(G,\bk)$ by left translation. For any natural transformation $\phi \colon \omega \Rightarrow \omega$, the component $\phi_{\fun^r(G)}$ should commute with all morphisms in $\End_{\rep(G)}(\fun^r(G)) \simeq \bk[G]$. Thus $\phi_{\fun^r(G)}$ is a left $\bk[G]$-module map from $\fun^l(G)$ to itself. Similarly, we also have $\End_{\rep(G)}(\fun^l(G)) \simeq \bk[G]$ acting on $\Hom(G,\bk)$ by right translation.
\item Viewing the commutative algebra $\fun^r(G) \in \rep(G)$ as a braided lax monoidal functor $\ast \to \rep(G)$, then for any $\phi \in \Aut^{\EE_3}(\omega)$ the component $\phi_{\fun^r(G)} \colon \omega(\fun^r(G)) \to \omega(\fun^r(G))$ is an algebra homomorphism. It follows that $\phi \mapsto \phi_{\fun^r(G)}$ defines a group homomorphism $\Aut^{\EE_3}(\omega) \to \Aut^{\EE_2}_G(\fun^l(G))$, where the latter is the automorphism group of the regular algebra $\fun^l(G) \in \rep(G)$. One can also compute that the automorphism group $\Aut^{\EE_2}_G(\fun^l(G))$ is isomorphic to $G$, acting on $\fun^l(G)$ by right translation. We denote the composite map $\Aut^{\EE_3}(\omega) \to \Aut^{\EE_2}_G(\fun^l(G)) \simeq G$ by $\Psi$.
\eit
Clearly $\Psi \circ \Phi \colon G \to G$ is identity. To show that $\Phi \circ \Psi \colon \Aut^{\EE_3}(\omega) \to \Aut^{\EE_3}(\omega)$ is identity, we need to use the following fact: the natural inclusion
\[
\mathrm B \End_{\rep(G)}(\fun^r(G)) \to \rep(G)
\]
is a Karoubi completion. By the universal property of the Karoubi completion, a natural transformation $\phi \colon \omega \Rightarrow \omega$ is completely determined by the component $\phi_{\fun^r(G)} \colon \omega(\fun^r(G)) \to \omega(\fun^r(G))$.
\epf

\section{2-groups and 2-representations} \label{sec:2-group}

\subsection{2-groups and their homotopy groups}

In this subsection, we recall some basic facts about 2-groups.

\begin{defn}
A \emph{2-group} (also called a \emph{weak 2-group} or a \emph{Gr-category} \cite{Sin75}) is a monoidal category in which every object is invertible and every morphism is an isomorphism.
\end{defn}

In particular, by Lemma \ref{lem:invertible_object}, a 2-group is a rigid monoidal category.

\begin{rem}
For a 2-group $\CG$, if we fix a left dual $(g^L,b_g \colon \one \to g \otimes g^L,d_g \colon g^L \otimes g \to \one)$ for every object $g \in \CG$, it is clear that $(g^L,d_g^{-1},b_g^{-1})$ is a right dual of $g$. One can also check that the corresponding left dual functor and the right dual functor are equal as monoidal functors $\CG^\op \to \CG^\rev$. Therefore, we obtain a pivotal structure on $\CG$. It is easy to see that both the left and right quantum dimension of an object is the unit $e \in \End_\CG(\one)$. However, in general, the left and right quantum trace of an endomorphism may not coincide.
\end{rem}

\begin{rem} \label{rem:coherent_2-group}
A \emph{coherent 2-group} \cite{BL04} is a 2-group $\CG$ equipped with a choice of a left dual functor
\[
(-)^L \colon \CG^\op \to \CG^\rev .
\]
This `contravariant' left dual functor has an equivalent `covariant' version \cite{BL04}: $(-)^* \colon \CG \to \CG^\rev$ defined by
\[
\bigl( g \xrightarrow{a} h \bigr)^* \coloneqq \bigl( g^* \coloneqq g^L \xrightarrow{a^* \coloneqq (a^L)^{-1}} h^L \eqqcolon h^* \bigr) .
\]
In other words, the functor $(-)^*$ is the composition of the left dual functor and the functor $\CG \to \CG^\op$ that acts on objects as identity and sends each morphism to its inverse. We call this functor $(-)^*$ the \emph{inverse functor} of $\CG$.

By Lemma \ref{lem:invertible_object}, it is not hard to see that the notion of a 2-group is equivalent to that of a coherent 2-group. More precisely, the 2-category of 2-groups and the 2-category of coherent 2-groups (with the 1-morphisms and 2-morphisms being the monoidal functors and monoidal natural transformations, respectively) are equivalent \cite{BL04}. In the following, we always assume that a 2-group is a coherent 2-group, that is, equipped with an inverse functor.
\end{rem}

\begin{rem}
A \emph{2-groupoid} is a 2-category in which every 1-morphism is an equivalence and every 2-morphism is an isomorphism. Clearly the notion of a 2-group is equivalent to that of a 2-groupoid with only one object. Given a 2-group $\CG$, we use $\mathrm B \CG$ to denote its one-point delooping, i.e., the 2-category with only one object $\ast$ and $\Hom(\ast,\ast) = \CG$.
\end{rem}

\begin{expl} \label{expl:auto-equivalence_2-group}
Let $\SC$ be a 2-category. For every object $x \in \SC$, the 1-auto-equivalences on $x$ and 2-isomorphisms between them form a 2-group, denoted by $\Aut_\SC(x)$. This 2-group is called the \emph{auto-equivalence 2-group} of $x$.
\end{expl}


\begin{expl} \label{expl:skeletal_2-group}
Let $G$ be a group, $A$ be an abelian group equipped with a $G$-action $G \to \Aut(A)$ and $\alpha \in Z^3(G;A)$ be a 3-cocycle with coefficients in $A$. We denote the $G$-action on $A$ by
\begin{align*}
G \times A & \to A \\
(g,a) & \mapsto g \triangleright a .
\end{align*}
Then we define a 2-group $\CG(G,A,\alpha)$ as follows.
\bit
\item The set of objects is the underlying set of $G$.
\item For $g,h \in G$, the hom space between them is defined by
\[
\Hom(g,h) \coloneqq \begin{cases} A , & g = h , \\ \emptyset , & g \neq h .\end{cases}
\]
\item The composition of morphisms is given by the multiplication of $A$.
\item The tensor product functor $\otimes \colon \CG(G,A,\alpha) \times \CG(G,A,\alpha) \to \CG(G,A,\alpha)$ is defined by
\[
\bigl( g \xrightarrow{a} g \bigr) \otimes \bigl( h \xrightarrow{b} h \bigr) \coloneqq \bigl( gh \xrightarrow{a (g \triangleright b)} gh \bigr) , \quad g,h \in G , \, a,b \in A .
\]
\item The associator is defined by $\alpha_{g,h,k} = \alpha(g,h,k) \in A$.
\item The left and right unitors are defined by $\lambda_g = \alpha(e,e,g)^{-1}$ and $\rho_g = \alpha(g,e,e)$.
\eit
Note that the pentagon equation is precisely the cocycle condition for $\alpha$:
\[
(g \triangleright \alpha(h,k,l)) \alpha(g,hk,l) \alpha(g,h,k) = \alpha(g,h,kl) \alpha(gh,k,l) .
\]

Suppose $F \colon \CG(G,A,\alpha) \to \CG(H,B,\beta)$ is a monoidal functor between two 2-groups. Then by definition $F$ is determined by the following data.
\bit
\item The underlying functor of $F$ is equivalent to a group homomorphism $f_1 \colon G \to H$ and a group homomorphism $f_2 \colon A \to B$, also denoted by $\pi_1(F)$ and $\pi_2(F)$, respectively.
\item The monoidal structure of $F$ is given by a 2-cochain $\gamma \in C^2(G;B)$ (where the $G$-action on $B$ is induced by the $H$-action on $B$ and $f_1$) satisfying $\mathrm d \gamma \cdot f_1^*\beta = f_2 \circ \alpha$.
\eit
It follows that there exists a monoidal equivalence $\CC(G,A,\alpha) \simeq \CC(G,A,\alpha')$ whose underlying functor is identity if and only if $[\alpha] = [\alpha'] \in H^3(G;A)$.

There are obvious monoidal functors $\mathrm B A \to \CG(G,A,\alpha) \to G$. Thus we also denote $\CG(G,A,\alpha)$ by $G \rtimes_\alpha \mathrm B A$ because it is a 2-group extension of $G$ and $\mathrm B A$.
\end{expl}

\begin{rem}
Consider monoidal functors from a group $G$ to a 2-group $\CG(H,A,\alpha)$. By Example \ref{expl:skeletal_2-group}, a group homomorphism $f \colon G \to H$ can be extended to a monoidal functor $F \colon G \to \CG(H,A,\alpha)$ if and only if the cohomology class $[f^* \alpha] \in H^3(G;A)$ vanishes. Thus $[f^* \alpha]$ is the obstruction class of extending $f$ to a 2-group homomorphism. Moreover, when the obstruction class vanishes, all the monoidal functors extending $f$, modulo monoidal natural isomorphisms, form a torsor over $H^2(G;A)$. When the 2-group $\CG(H,A,\alpha)$ is the braided autoequivalence 2-group of a nondegenerate braided fusion category, this obstruction class is exactly the $\mathcal O_3$ obstruction in \cite{ENO10}.

On the other hand, we always have a general type of extension. Let $\omega \coloneqq f^* \alpha \in Z^3(G;A)$. Then there exists a monoidal functor $F \colon \CG(G,A,\omega) = G \rtimes_\omega A \to \CG(H,A,\alpha)$ extending $f$, in the sense that $\pi_1(F) \colon G \to H$ is $f$ and $\pi_2(F) \colon A \to A$ is identity. Moreover, all such monoidal functors, modulo monoidal natural isomorphisms, form a torsor over $H^2(G;A)$.
\end{rem}

\begin{defn}
Let $\CG$ be a 2-group. Its \emph{first homotopy group} $\pi_1(\CG)$ is the group of isomorphism classes of objects of $\CG$ with the multiplication induced by the tensor product of $\CG$. Its \emph{second homotopy group} $\pi_2(\CG)$ is the endomorphism group $\CG(\one,\one)$ of the tensor unit $\one$.
\end{defn}

By the Eckmann-Hilton argument, the second homotopy group $\pi_2(\CG)$ of any 2-group $\CG$ is always an abelian group.



Let $\CG$ be a 2-group. Recall that a category is \emph{skeletal} if every isomorphism class of objects contains only one object. Since every category is equivalent to a skeletal one, we can find a skeletal category that is equivalent to $\CG$ as categories and transfer the 2-group structure of $\CG$ to the skeletal one. This argument shows that every 2-group is equivalent to a skeletal 2-group. 
Now we assume that $\CG$ is skeletal. Then $\pi_1(\CG)$ naturally acts on $\pi_2(\CG)$ by
\[
g \triangleright a \coloneqq (1_g \otimes a) \otimes 1_{g^{-1}} , \quad g \in \pi_1(\CG) , \, a \in \pi_2(\CG) = \CG(e,e) .
\]
For every $g,h,k \in \pi_1(\CG)$, we define
\[
\alpha(g,h,k) \coloneqq \alpha_{g,h,k} \otimes 1_{(ghk)^{-1}} \in \CG(e,e) = \pi_2(\CG) .
\]
Then the pentagon equation reads
\[
(g \triangleright \alpha(h,k,l)) \alpha(g,hk,l) \alpha(g,h,k) = \alpha(g,h,kl) \alpha(gh,k,l) .
\]
In other words, $\alpha \in Z^3(\pi_1(\CG);\pi_2(\CG))$ is a 3-cocycle.

The homotopy groups $\pi_1(\CG),\pi_2(\CG)$, together with the $\pi_1(\CG)$-action on $\pi_2(\CG)$ and the 3-cohomology class $[\alpha] \in H^3(\pi_1(\CG);\pi_2(\CG))$, are invariant under equivalences of 2-groups. The definition of these data for non-skeletal 2-groups can be found in \cite{Sin75}.

Conversely, from the data $(\pi_1(\CG),\pi_2(\CG),\alpha)$ we can also construct a 2-group $\CG(\pi_1(\CG),\pi_2(\CG),\alpha)$ (see Example \ref{expl:skeletal_2-group}). By the above discussion, it is not hard to see that $\CG$ is monoidally equivalent to $\CG(\pi_1(\CG),\pi_2(\CG),\alpha)$. So we have the following classification result.

\begin{thm}[\cite{Sin75}]
A 2-group $\CG$ is completely determined by its homotopy groups $\pi_1(\CG),\pi_2(\CG)$, together with the conjugation $\pi_1(\CG)$-action on $\pi_2(\CG)$ and the 3-cohomology class $[\alpha] \in H^3(\pi_1(\CG);\pi_2(\CG))$ represented by the associator.
\end{thm}

\begin{defn}
A \emph{finite 2-group} is a 2-group $\CG$ such that both $\pi_1(\CG)$ and $\pi_2(\CG)$ are finite groups.
\end{defn}

%

\subsection{2-representations of finite 2-groups}

We use $2\vect$ to denote the 2-category of finite semisimple ($\bk$-linear) categories.

\begin{defn}
Let $\CG$ be a 2-group. A \emph{finite semisimple 2-representation} of $\CG$ is a 2-functor $\mathrm B \CG \to 2\vect$, where $\mathrm B \CG$ is the one-point delooping 2-category of $\CG$.
\end{defn}

\begin{rem}
Given a 2-group $\CG$, the following notions are equivalent:
\bit
\item a 2-functor $\mathrm B \CG \to 2\vect$;
\item a finite semisimple category $\CV$ equipped with a monoidal functor $\CG \to \End_\bk(\CV)$, where $\End_\bk(\CV)$ is the monoidal category of $\bk$-linear endofunctors on $\CV$ and natural transformations;
\item a finite semisimple category $\CV$ equipped with a monoidal functor $\CG \to \Aut_\bk(\CV)$, where $\Aut_\bk(\CV)$ is the 2-group of $\bk$-linear auto-equivalence on $\CV$ and natural isomorphisms;
\item a finite semisimple left $\CG$-module, i.e., a finite semisimple category $\CV$ equipped with a left $\CG$-action $\CG \times \CV \to \CV$ and coherence data.
\eit
More generally, one can similarly define the notion of a 2-representation of $\CG$ in any 2-category \cite{Elg07}.
\end{rem}

\begin{expl}
The finite semisimple category $\vect$ equipped with the trivial $\CG$-action is called the \emph{trivial 2-representation} of $\CG$. In other words, the $\CG$-action is given by the composite functor $\CG \to \{\ast\} \to \End_\bk(\CV)$ where $\{\ast\}$ is the trivial monoidal category with only one object and identity morphism.
\end{expl}

\begin{expl} \label{expl:regular_2-rep}
Suppose $\CG$ is a finite 2-group. The functor category $\fun(\CG,\vect)$ has a left $\CG$-action defined by the right translation:
\[
(g \odot^r F)(h) \coloneqq F(h \otimes g) , \quad g,h \in \CG , \, F \in \fun(\CG,\vect) .
\]
Since $\CG$ is equivalent to the disjoint union of $\lvert \pi_1(\CG) \rvert$ copies of $\mathrm B \pi_2(\CG)$, the functor category $\fun(\CG,\vect)$ is equivalent to the direct sum of $\lvert \pi_1(\CG) \rvert$ copies of $\rep(\pi_2(\CG))$ and hence finite semisimple. We denote this finite semisimple 2-representation by $\mathcal F_r(\CG) \in 2\rep(\CG)$. There is another left $\CG$-action on $\fun(\CG,\vect)$ defined by the left translation:
\[
(g \odot^l F)(h) \coloneqq F(g^* \otimes h) , \quad g,h \in \CG , \, F \in \fun(\CG,\vect) ,
\]
where $(-)^* \colon \CG \to \CG^\rev$ is the inverse functor (see Remark \ref{rem:coherent_2-group}). We denote this 2-representation by $\mathcal F_l(\CG) \in 2\rep(\CG)$. These two 2-representations are equivalent because the left translation on $\fun(\CG,\vect)$ is the same as the right translation on $\fun(\CG^\rev,\vect)$ via the inverse functor. We call them the \emph{left and right regular 2-representations} of $\CG$, respectively (see also \cite{Elg11}).
\end{expl}

\begin{expl} \label{expl:dual_2-representation}
Let $\CV \in 2\rep(\CG)$. The functor category $\fun_\bk(\CV,\vect)$ is naturally a right $\CG$-module:
\[
(F \odot g)(v) \coloneqq F(g \odot v) , \quad F \in \fun_\bk(\CV,\vect) , \, g \in \CG , \, v \in \CV .
\]
Therefore, $\fun_\bk(\CV,\vect)$ is a left $\CG^\rev$-module and can be pulled back along the inverse functor $(-)^* \colon \CG \to \CG^\rev$ (see Remark \ref{rem:coherent_2-group}) to be a left $\CG$-module:
\[
(g \odot F)(v) \coloneqq F(g^* \odot v) , \quad F \in \fun_\bk(\CV,\vect) , \, g \in \CG , \, v \in \CV .
\]
Since $\CV$ is finite semisimple, the Yoneda embedding $\CV^\op \to \fun(\CV,\vect)$ is an equivalence. So the left $\CG$-module structure on $\fun_\bk(\CV,\vect)$ can be transferred to $\CV^\op$ along this equivalence. Equivalently, $\CV^\op$ is a left $\CG^\op$-module and thus can be pulled back along the inverse morphism functor $\CG \to \CG^\op$ (see Remark \ref{rem:coherent_2-group}) to be a left $\CG$-module. We call $\fun_\bk(\CV,\vect)$ or $\CV^\op$ the \emph{dual 2-representation} of $\CV$.
\end{expl}


\begin{defn}
The \emph{2-category of finite semisimple 2-representations} of $\CG$ is the 2-category $\fun(\mathrm B\CG,2\vect)$. We denote this 2-category by $2\rep(\CG)$.
\end{defn}

By viewing the objects in $2\rep(\CG)$ as finite semisimple left $\CG$-modules, the 1-morphisms and 2-morphisms in $2\rep(\CG)$ are left $\CG$-module functors and left $\CG$-module natural transformations, respectively.


\subsection{2-group 2-algebra}

%
%
%

Recall that for a group $G$, the notion of a $G$-representation is equivalent to that of a module over the group algebra $\bk[G]$. In this subsection we categorify this idea and define a `2-group 2-algebra', denoted by $\vect_\CG$, for a finite 2-group $\CG$.

Let $\CG$ be a finite 2-group. We define its \emph{2-group 2-algebra} $\vect_\CG$ in two steps. First, recall that for any category $\CC$ we can linearize its hom spaces and obtain a $\bk$-linear category $\bk_\ast \CC$ as follows:
\bit
\item $\ob(\bk_\ast \CC) \coloneqq \ob(\CC)$;
\item $\bk_\ast \CC(x,y) \coloneqq \bk[\CC(x,y)]$ for any $x,y \in \CC$, where $\bk[S]$ denotes the free vector space generated by a set $S$;
\item the composition and identities of $\bk_\ast \CC$ are induced by those of $\CC$.
\eit
Also, there is a canonical inclusion functor $\CC \to \bk_\ast \CC$. Moreover, if $\CC$ is a monoidal category, then $\bk_\ast \CC$ is also a $\bk$-linear monoidal category and the inclusion functor $\CC \to \bk_\ast \CC$ is monoidal.

\begin{rem}
In general, every lax monoidal functor $F \colon \CA \to \CB$ between two monoidal categories induces a 2-functor from the 2-category of $\CA$-enriched categories to that of $\CB$-enriched categories. Here we use the free vector space functor $\bk[-] \colon \set \to \vect$ to construct $\bk_\ast \CC$.
\end{rem}

Then we define $\vect_\CG \coloneqq \Kar(\bk_\ast \CG)$ to be the Karoubi completion of $\bk_\ast \CG$. Since $\bk_\ast \CG$ is equivalent to the the disjoint union of $\lvert \pi_1(\CG) \rvert$ copies of $\mathrm B \bk[\pi_2(\CG)]$, by Example \ref{expl:Karoubi_completion_group} the Karoubi completion $\Kar(\bk_\ast \CG)$ is equivalent to the direct sum of $\lvert \pi_1(\CG) \rvert$ copies of $\rep(\pi_2(\CG))$ and hence finite semisimple.

\begin{rem} \label{rem:vect_G_dual}
When $\CG = G$ is a finite 1-group, $\vect_G$ is usually defined as the category of finite-dimensional $G$-graded vector spaces, or the functor category $\fun(G,\vect)$. Thus an alternative definition of $\vect_\CG$ is the functor category $\fun(\CG,\vect)$, in which the objects can be viewed as `$\CG$-graded vector spaces'. By the universal property of the free vector space functor $\bk[-] \colon \set \to \vect$, the functor category $\fun(\CG,\vect)$ is equivalent to the $\bk$-linear functor category $\fun_\bk(\bk_\ast \CG,\vect)$. Then by the universal property of the Karoubi completion, this $\bk$-linear functor category is also equivalent to $\fun_\bk(\Kar(\bk_\ast \CG),\vect)$. Recall for any semisimple category $\CC$ the Yoneda embedding $\CC^\op \to \fun_\bk(\CC,\vect)$ defined by $x \mapsto \CC(x,-)$ is a $\bk$-linear equivalence. So $\fun(\CG,\vect)$ is equivalent to $\Kar(\bk_\ast \CG)^\op$ as $\bk$-linear categories. Thanks to the inverse morphism functor $\CG \to \CG^\op$ (see Remark \ref{rem:coherent_2-group}), $\Kar(\bk_\ast \CG)^\op$ is equivalent to $\Kar(\bk_\ast \CG)$. The composite equivalence
\[
\vect_\CG \simeq \vect_\CG^\op \simeq \fun_\bk(\vect_\CG,\vect) \simeq \fun(\CG,\vect)
\]
categorifies the Fourier transform of finite group (see \cite[Section 3.4]{HXZ24} for more details). Here we choose to define $\vect_\CG$ to be $\Kar(\bk_\ast \CG)$ because there is a natural inclusion $\CG \to \Kar(\bk_\ast \CG)$, but the corresponding functor $\CG \to \fun(\CG,\vect)$ is not so natural.
\end{rem}

The monoidal structure of $\CG$ naturally induces that of $\Kar(\bk_\ast \CG) = \vect_\CG$ so that $\vect_\CG$ is a multi-fusion category and the natural inclusion functor $\CG \to \vect_\CG$ is monoidal (see Section \ref{sec_Kar}). Let us list the simple objects and fusion rules of $\vect_\CG$. Assume that $\CG = \CG(G,A,\alpha)$ is skeletal. Then the simple objects of $\vect_\CG$ are pairs $(g,p_\rho)$, where $g \in G$, $\rho \in \hat A$ and $p_\rho \in \bk[A]$ is the idempotent defined by
\[
p_\rho \coloneqq \frac{1}{\lvert A \rvert} \sum_{a \in A} \rho(a)^{-1} \cdot a .
\]
The vector $p_\rho$ spans a subrepresentation of $\bk[A]$ which is isomorphic to $\rho$. For $\rho,\sigma \in \hat A$, we have $p_\rho p_\sigma = \delta_{\rho,\sigma} p_\sigma$. So these idempotents induce an isomorphism between the group algebra $\bk[A]$ and the function algebra $\fun(\hat A)$ (this is also the Fourier transform). Moreover, the $G$-action on $A$ transfers to $\fun(\hat A)$ and hence $\hat A$ along this isomorphism:
\[
(g \triangleright \rho)(a) \coloneqq \rho(g^{-1} \triangleright a) , \quad g \in G , \, \rho \in \hat A , \, a \in A .
\]
Then one can compute that the fusion rules are given by the following formula:
\[
(g,p_\rho) \otimes (h,p_\sigma) = (gh,p_\rho(g \triangleright p_\sigma)) = (gh,p_\rho p_{g \triangleright \sigma}) = \delta_{\rho,g \triangleright \sigma} (gh,p_\rho) , \quad g,h \in G , \, \rho,\sigma \in \hat A .
\]
For simplicity, we denote the simple object $(g,p_\rho)$ by $(g,\rho)$ if there is no ambiguity. By \ref{lem:idempotent_completion_monoidal}, the associator
\[
\alpha_{(g,\rho),(h,\sigma),(k,\tau)} \colon ((g,\rho) \otimes (h,\sigma)) \otimes (k,\tau) \to (g,\rho) \otimes ((h,\sigma) \otimes (k,\tau))
\]
is given by the composite morphism
\[
(p_\rho \otimes (p_\sigma \otimes p_\tau)) \circ \alpha(g,h,k) \circ ((p_\rho \otimes p_\sigma) \otimes p_\tau) = p_\rho p_{g \triangleright \sigma} p_{gh \triangleright \tau} \cdot \alpha(g,h,k) .
\]
The tensor unit is the direct sum
\[
\one = \bigoplus_{\rho \in \hat A} (e,\rho) .
\]
The left and right dual of $(g,\rho)$ is $(g^{-1},g^{-1} \triangleright \rho)$.

Now we apply Theorem \ref{thm:structure_multi-fusion_category} to $\vect_\CG$. The simple objects in $(\vect_\CG)_{\rho,\sigma} \coloneqq (e,\rho) \otimes \vect_\CG \otimes (e,\sigma)$ are
\[
\{(g,\rho) \mid g \in G , \, g \triangleright \sigma = \rho\} .
\]
Denote the stabilizer subgroup of $\rho \in \hat A$ by $\Stab(\rho) \subseteq G$. Then $(\vect_\CG)_{\rho,\rho}$ is a pointed fusion category whose group of (isomorphism classes of) simple objects is $\Stab(\rho)$, and the associator is given by
\[
\alpha_{(g,\rho),(h,\rho),(k,\rho)} = p_\rho \cdot \alpha(g,h,k) = \rho(\alpha(g,h,k)) \cdot p_\rho .
\]
Note that $p_\rho$ is the identity morphism on $(g,\rho)$. Thus the 3-cocycle representing the associator of the pointed fusion category $(\vect_\CG)_{\rho,\rho}$ is $(\rho \circ \alpha) \vert_{\Stab(\rho)} \in Z^3(\Stab(\rho);\bk^\times)$.

\begin{prop} \label{prop_vect_G_structure}
Let $\CG$ be a finite 2-group. Then $\vect_\CG$ is a multi-fusion category, and the indecomposable direct summands of $\vect_\CG$ are one-to-one corresponding to the $\pi_1(\CG)$-orbits in $\widehat{\pi_2(\CG)}$. For $\rho \in \widehat{\pi_2(\CG)}$, the indecomposable direct summand corresponding to $[\rho]$ is Morita equivalent to $\vect_{\Stab(\rho)}^{(\rho \circ \alpha)\vert_{\Stab(\rho)}}$.
\end{prop}

\begin{rem} \label{rem:equivalence_regular_2-rep_dual}
The tensor product of $\CG$ induces a left $\CG$-module structure on $\vect_\CG$, i.e., $\vect_\CG \in 2\rep(\CG)$. The dual 2-representation of $\vect_\CG$ is the regular 2-representation $\mathcal F_l(\CG)$ defined in Example \ref{expl:regular_2-rep}. It follows from Remark \ref{rem:vect_G_dual} that $\vect_\CG$ is also equivalent to the regular 2-representation $\mathcal F_l(\CG)$ or $\mathcal F_r(\CG)$.
\end{rem}

\begin{prop} \label{prop:2-rep_module_over_2-group_algebra}
There is an equivalence of $\bk$-linear 2-categories $2\rep(\CG) \simeq \LMod_{\vect_\CG}(2\vect)$.
\end{prop}

\pf
By the universal property of the free vector space functor and that of the Karoubi completion, we have
\[
2\rep(\CG) \simeq \fun(\mathrm B \CG,2\vect) \simeq \fun_\bk(\mathrm B \bk_\ast \CG,2\vect) \simeq \fun_\bk(\mathrm B \Kar(\bk_\ast \CG),2\vect) \simeq \LMod_{\Kar(\bk_\ast \CG)}(2\vect) . \qed
\]
\epf

\begin{rem}
Over a general ground field $\bk$, the category $\vect_\CG$ may not be semisimple. But the notion of a finite semisimple $\CG$-module is still equivalent to that of a finite semisimple $\vect_\CG$-module, even though $\vect_\CG$ is not an object in $2\vect$.
\end{rem}

\subsection{The 2-category of 2-representations} \label{sec:2Rep(G)_2-cat}

Let $\CG$ be a finite 2-group. By Theorem \ref{thm:semisimple_2-cat_multi-fusion}, the 2-category of finite semisimple modules over a multi-fusion category is a finite semisimple 2-category. Hence from Proposition \ref{prop:2-rep_module_over_2-group_algebra} we obtain the following result.

\begin{cor} \label{cor:2-representation_semisimple}
Let $\CG$ be a finite 2-group. Then $2\rep(\CG)$ is a finite semisimple 2-category.
\end{cor}

The above result was sketchily proved in \cite[Example 1.4.19]{DR18}. See also \cite[Theorem 5.8]{KZ24} for higher-categorical generalizations. Also by Proposition \ref{prop_vect_G_structure} we immediately obtain the following result.

\begin{thm} \label{thm_2Rep_G_structure}
Let $\CG$ be a finite 2-group. Let $G \coloneqq \pi_1(\CG)$, $A \coloneqq \pi_2(\CG)$ and $\alpha \in Z^3(G;A)$ be the 3-cocycle representing the associator of $\CG$. Then we have the equivalence of finite semisimple 2-categories
\[
2\rep(\CG) \simeq \bigoplus_{[\rho] \in \hat A/G} 2\rep(\Stab(\rho),(\rho \circ \alpha)\vert_{\Stab(\rho)}) .
\]
In particular, the connected components of $2\rep(\CG)$ are one-to-one corresponding to $\pi_1(\CG)$-orbits in $\widehat{\pi_2(\CG)}$.
\end{thm}

\begin{expl}
Let $\CG$ be a finite 2-group. If $\pi_2(\CG)$ is trivial, then $\vect_\CG \simeq \vect_{\pi_1(\CG)}$ as fusion categories, and thus $2\rep(\CG)$ is equivalent to $2\rep(\pi_1(\CG))$ as finite semisimple 2-categories. If $\pi_1(\CG)$ is trivial, then $\vect_\CG \simeq \vect^{\oplus \widehat{\pi_2(\CG)}}$ as multi-fusion categories, and thus $2\rep(\CG)$ is equivalent to $2\vect_{\widehat{\pi_2(\CG)}}$ as finite semisimple 2-categories.
\end{expl}

\begin{expl} \label{expl:trivial_component_2Rep}
Let $\CG = \CG(G,A,\alpha)$ be a skeletal 2-group. Note that the $G$-orbit of the trivial character $1 \in \hat A$ contains only one element. So the fusion category $(\vect_\CG)_{1,1}$ is a direct summand of $\vect_\CG$. By the above discussion, $(\vect_\CG)_{1,1}$ is equivalent to $\vect_{\pi_1(\CG)}$ as fusion categories. It follows that $2\rep(\CG)$ always has a connected component equivalent to $2\rep(\pi_1(\CG))$. We call it the \emph{trivial connected component} of $2\rep(\CG)$. The trivial 2-representation $\vect$ is a simple object in this connected component, and the endomorphism category of the trivial 2-representation is $\End_{2\rep(\CG)}(\vect) \simeq \fun_{\vect_{\pi_1(\CG)}}(\vect,\vect) \simeq \rep(\pi_1(\CG))$.
\end{expl}

\begin{expl} \label{expl:endomorphism_regular_2-rep}
The endomorphism category of $\vect_\CG \in 2\rep(\CG)$ is $\End_{2\rep(\CG)}(\vect_\CG) = \fun_{\vect_\CG}(\vect_\CG,\vect_\CG) \simeq \vect_\CG^\rev$. By carefully tracing out the equivalence $\vect_\CG \simeq \mathcal F_l(\CG)$ (see Remark \ref{rem:equivalence_regular_2-rep_dual}), it is easy to see that the endomorphism category $\End_{2\rep(\CG)}(\mathcal F_l(\CG))$ is equivalent to $\vect_\CG$, which acts on $\mathcal F_l(\CG)$ by right translation. Similarly, $\End_{2\rep(\CG)}(\mathcal F_r(\CG)) \simeq \vect_\CG$ acts on $\mathcal F_r(\CG)$ by left translation.
\end{expl}

\begin{defn}
Let $\CG$ be a finite 2-group. The \emph{tensor product} of two finite semisimple 2-representations $\CV,\CW \in 2\rep(\CG)$ is the Deligne tensor product $\CV \boxtimes \CW$ equipped with the $\CG$-action defined by
\[
g \odot (v \boxtimes w) \coloneqq (g \odot v) \boxtimes (g \odot w) , \quad g \in \CG , \, v \in \CV , \, w \in \CW .
\]
\end{defn}

The tensor product of 2-representations equips $2\rep(\CG)$ a $\bk$-linear symmetric monoidal 2-category structure, which is induced by the universal property of the Deligne tensor product of finite semisimple categories.

\begin{thm}
Let $\CG$ be a finite 2-group. Then $2\rep(\CG)$ is a symmetric fusion 2-category.
\end{thm}

\pf
By Corollary \ref{cor:2-representation_semisimple}, $2\rep(\CG)$ is a finite semisimple symmetric monoidal 2-category. Clearly the tensor unit, i.e., the trivial 2-representation of $\CG$, is simple. It remains to prove that $2\rep(\CG)$ admits left and right duals for objects.

Suppose $\CV \in 2\rep(\CG)$. Let us show that the dual 2-representation defined in Example \ref{expl:dual_2-representation} is a left dual of $\CV$. Taking the dual 2-representation to be $\CV^\op$, the evaluation 1-morphism is the hom functor $\Hom \colon \CV^\op \boxtimes \CV \to \vect$, and the coevaluation 1-morphism $\vect \to \CV \boxtimes \CV^\op$ is defined by $\bk \mapsto \bigoplus_{x \in \Irr(\CV)} x \boxtimes x$. Taking the dual 2-representation to be $\fun_\bk(\CV,\vect)$, the evaluation 1-morphism is the evaluation functor, and the coevaluation 1-morphism $\vect \to \CV \boxtimes \fun_\bk(\CV,\vect)$ is defined by $\bk \mapsto \bigoplus_{x \in \Irr(\CV)} x \boxtimes \CV(x,-)$. It is easy to verify that these functors are indeed 1-morphisms in $2\rep(\CG)$ and the zig-zag equations hold. The proof of that the dual 2-representations are right duals is similar.
\epf

\begin{rem}
Since $\CV$ is finite semisimple over an algebraically closed field, the coevaluation 1-morphism can be equivalently defined by using (co)ends: $\bk \mapsto \int_{x \in \CV} x \boxtimes x$ or $\bk \mapsto \int^{x \in \CV} x \boxtimes x$.
\end{rem}

For $\CV,\CW \in 2\rep(\CG)$, their Deligne tensor product $\CV \boxtimes \CW$ is naturally a left $\CG \times \CG$-module and can be pulled back along the diagonal monoidal functor $\CG \to \CG \times \CG$ to be a left $\CG$-module. This gives an equivalent definition of the tensor product of $2\rep(\CG)$. Indeed, the 2-group $\CG$ equipped with the diagonal map $\CG \to \CG \times \CG$, the constant functor $\CG \to \ast$ and the inverse functor $\CG^\rev \to \CG$ is a (categorified) Hopf algebra object in the 2-category of categories. Then $\vect_\CG$, as the `linearization' of $\CG$, is a (finite semisimple) \emph{Hopf monoidal category} (see for example \cite{CF94,Neu97}). The Hopf monoidal structure on $\vect_\CG$ naturally induces a monoidal structure on the 2-category of finite semisimple $\vect_\CG$-modules, which is the same as the monoidal structure of $2\rep(\CG)$.

Moreover, by the Tannaka-Krein reconstruction for fusion 2-categories \cite{Gre23}, the endomorphism category of the forgetful 2-functor $\omega \colon 2\rep(\CG) \to 2\vect$ admits a canonical Hopf monoidal structure such that $2\rep(\CG) \simeq 2\rep(\End(\omega))$ as fusion 2-categories. We claim that $\End(\omega)$ is equivalent to $2\vect_\CG$. Indeed, for every $\CV \in 2\rep(\CG)$ and $v \in \CV$, there is a unique (up to equivalence) left $\CG$-module functor $\CG \to \CV$ sending the tensor unit $e \in \CG$ to $v$, that is, the functor $- \odot v \colon \CG \to \CV$. Thus there is a canonical equivalence
\[
2\rep(\CG)(\vect_\CG,\CV) \simeq \omega(\CV) \simeq 2\vect(\vect,\omega(\CV)) .
\]
Therefore, we have $\vect_\CG \simeq \omega^L(\vect)$, where $\omega^L \colon 2\vect \to 2\rep(\CG)$ is the left adjoint of $\omega$. Also, $\vect_\CG \in 2\rep(\CG)$ represents the forgetful 2-functor $\omega$ (this was also obtained by Elgueta \cite{Elg11}). By the Yoneda lemma, we have a monoidal equivalence
\[
\End(\omega) \simeq \End(2\rep(\CG)(\vect_\CG,-)) \simeq 2\rep(\CG)(\vect_\CG,\vect_\CG)^\rev \simeq \vect_\CG .
\]
This equivalence sends $g \in \CG \hookrightarrow \vect_\CG$ to the 2-natural transformation $\Phi^g \colon \omega \Rightarrow \omega$ defined by $\Phi^g_\CV \coloneqq (g \odot -) \colon \omega(\CV) \to \omega(\CV)$. The cotensor product $\Delta$ on $\End(\omega)$ gives
\[
\Delta(\Phi^g)_{\CV \boxtimes \CW} \coloneqq \Phi^g_{\CV \boxtimes \CW} = (g \odot -) \boxtimes (g \odot -) \colon \omega(\CV) \boxtimes \omega(\CW) \to \omega(\CV) \boxtimes \omega(\CW) ,
\]
where we identify $\End(\omega) \boxtimes \End(\omega)$ with $\End(\omega \boxtimes \omega)$, and the first $\CV \boxtimes \CW$ is an object in $2\rep(\CG) \boxtimes 2\rep(\CG)$ while the second one is in $2\rep(\CG)$. Therefore, we see that the comonoidal structure on $\End(\omega)$ coincides with that on $\vect_\CG$. Similarly, by considering the dual in $2\rep(\CG)$, it is easy to see that the antipode on $\End(\omega)$ coincides with that on $\vect_\CG$.

\begin{rem} \label{rem:diagonal_action_simple}
Assume that $\CG = \CC(G,A,\alpha)$ is skeletal. 
The image of $(g,\rho) \in \vect_\CG$ under the cotensor product $\Delta \colon \vect_\CG \to \vect_\CG \boxtimes \vect_\CG$ is the retract of the idempotent
\[
\frac{1}{\lvert A \rvert} \sum_{a \in A} \rho(a)^{-1} \cdot a \otimes a \in (\vect_\CG \boxtimes \vect_\CG)(g \boxtimes g) \simeq  \bk[A] \otimes \bk[A] .
\]
A direct computation shows that this idempotent is equal to $\sum_{\sigma \tau = \rho} p_\sigma \otimes p_\tau$. Hence
\[
\Delta(g,\rho) \simeq \bigoplus_{\sigma \tau = \rho} (g,\sigma) \boxtimes (g,\tau) \in \vect_\CG \boxtimes \vect_\CG .
\]
This also gives the action of a simple object $(g,\rho) \in \vect_\CG$ on the tensor product of two 2-representations.
\end{rem}

\begin{expl}
Let $\CG$ be a finite 2-group. If $\pi_2(\CG)$ is trivial, $2\rep(\CG)$ is equivalent to $2\rep(\pi_1(\CG))$ as symmetric fusion 2-categories. If $\pi_1(\CG)$ is trivial, $2\rep(\CG)$ is equivalent to $2\vect_{\widehat{\pi_2(\CG)}}$ as symmetric fusion 2-categories.
\end{expl}

\begin{expl}
Let $\CG$ be a finite 2-group. Recall Example \ref{expl:trivial_component_2Rep} that the trivial connected component of $2\rep(\CG)$ is equivalent to $2\rep(\pi_1(\CG))$ as finite semisimple 2-categories. By Remark \ref{rem:diagonal_action_simple}, it is not hard to see that the trivial connected component is closed under the tensor product and equivalent to $2\rep(\pi_1(\CG))$ as symmetric fusion 2-categories.
\end{expl}

\begin{expl}
Let $\CG$ be a finite 2-group. Denote $G \coloneqq \pi_1(\CG)$ and $A \coloneqq \pi_2(\CG)$. Assume that the $G$-action on $A$ is trivial and the 3-cocycle representing the associator is trivial. Then each $G$-orbit in $\hat A$ consists of only one element and thus $\vect_\CG \simeq (\vect_G)^{\oplus \hat A}$ as multi-fusion categories. Therefore, $2\rep(\CG)$ is equivalent to $2\rep(G)^{\oplus \hat A}$ as finite semisimple 2-categories. In particular, $\pi_0(2\rep(\CG)) \simeq \hat A$,and each connected component is equivalent to $2\rep(G)$.

Moreover, in this case $\CG$ is equivalent to the Cartesian product $G \times \mathrm B A$ as 2-groups. Then $\vect_\CG \simeq \vect_G \boxtimes \vect_{\mathrm B A} \simeq \vect_G \boxtimes \vect^{\oplus \hat A}$ as multi-fusion categories. Thus $2\rep(\CG)$ is equivalent to $2\rep(G) \boxtimes 2\vect_{\hat A}$ as fusion 2-categories.
\end{expl}

\begin{expl} \label{expl_Z2_Z3_2-rep}
Let $\CG = \CG(\Zb_2,\Zb_3,1)$ be the finite 2-group defined as follows (recall Example \ref{expl:skeletal_2-group}):
\bit
\item Its first homotopy group is $\Zb_2$, and the second homotopy group is $\Zb_3$.
\item The action of $\Zb_2$ on $\Zb_3$ is the nontrivial one.
\item The 3-cocycle representing the associator is trivial. Indeed, in this case we have $H^3(\Zb_2;\Zb_3) = 0$.
\eit
Let us denote $\pi_1(\CG) = \{e,x\}$ with $x^2 = e$ and $\pi_2(\CG) = \{1,y,y^2\}$ with $y^3 = 1$. The elements in $\widehat{\pi_2(\CG)}$ are denoted by $\rho_k$ for $k = 0,1,2$, where $\rho_k(y) = \omega^k$ and $\omega \in \bk^\times$ is a 3rd root of unity. The $\pi_1(\CG)$-action on $\widehat{\pi_2(\CG)}$ is given by $x \triangleright \rho_1 = \rho_2$ and $x \triangleright \rho_2 = \rho_1$.

There are 6 simple objects in $\vect_\CG$. The simple objects $(e,\rho_0)$ and $(x,\rho_0)$ generate a fusion subcategory equivalent to $\vect_{\Zb_2}$, and the simple objects $(e,\rho_1),(x,\rho_1),(x,\rho_2),(e,\rho_2)$ generate a multi-fusion subcategory equivalent to $\End_\bk(\vect \oplus \vect)$. Therefore, we have $\vect_\CG \simeq \vect_{\Zb_2} \oplus \End_\bk(\vect \oplus \vect)$ as multi-fusion categories and $2\rep(\CG) \simeq 2\rep(\Zb_2) \oplus 2\vect$ as finite semisimple 2-categories. In particular, there are 3 simple 2-representations of $\CG$.
\bit
\item The first two are two indecomposable finite semisimple $\vect_{\Zb_2}$-modules. One is $\vect$, which is also the trivial 2-representation of $\CG$. We denote it by $\one$.
\item Another one is the regular $\vect_{\Zb_2}$-module $\vect_{\Zb_2}$. The object $x \in \CG$ acts on $\vect_{\Zb_2}$ by permuting two simple objects, and the morphism $y \in \pi_2(\CG)$ acts on $\vect_{\Zb_2}$ as the identity natural transformation. We denote this 2-representation by $\one_c$.
\item The last one is the unique indecomposable finite semisimple $\End_\bk(\vect \oplus \vect)$-module $\vect \oplus \vect$. The object $x \in \CG$ acts on $\vect \oplus \vect$ by permuting two simple objects, and the morphism $y \in \pi_2(\CG)$ acts on $\vect \oplus \vect$ as the natural transformation whose components on two simple objects are $\rho_1(y) = \omega$ and $\rho_2(y) = \omega^2$. There are two choices of this natural transformation, but they lead to equivalent 2-representations of $\CG$; the equivalence is given by permuting two simple objects. We denote this 2-representation by $S$.
\eit
This 2-category can be depicted by the following graph:
\[
\xymatrix{
\one \ar@(ul,ur)[]^{\rep(\Zb_2)}  \ar@/^/[rr]^{\vect} & & \one_c \ar@(ul,ur)[]^{\vect_{\Zb_2}} \ar@/^/[ll]^{\vect} & S \ar@(ul,ur)[]^{\vect}
} .
\]
Let us compute the fusion rule. For example, $\one_c \boxtimes S$ is a finite semisimple category with 4 simple objects, denoted by $e \boxtimes e,e \boxtimes x,x \boxtimes e,x \boxtimes x$. By considering the action of $x \in \CG$, we see that $e \boxtimes e$ and $x \boxtimes x$ generate a submodule, and $e \boxtimes x$ and $x \boxtimes e$ generate another submodule. Then by considering the action of $y$ we see that both two submodules are equivalent to $S$. Thus $\one_c \boxtimes S \simeq S \oplus S$. Similarly we obtain the fusion rules of $2\rep(\CG)$:
\[
\one_c \boxtimes \one_c \simeq \one_c \oplus \one_c , \quad S \boxtimes S \simeq \one_c \oplus S , \quad \one_c \boxtimes S \simeq S \boxtimes \one_c \simeq S \oplus S .
\]
This 2-group $\CG$ also appeared in Delcamp's work \cite{Del22a} as the input data of a 3+1D lattice (higher) gauge theory. Moreover, Delcamp showed that the topological defects on a certain boundary form a fusion 2-category which is equivalent to $2\rep(\CG)$.\footnote{Delcamp found $4$ simple objects $\one,\mathcal E_1,\mathcal E_\omega,\mathcal E_{\bar \omega}$, but $\mathcal E_\omega$ and $\mathcal E_{\bar \omega}$ are indeed equivalent. He checked this equivalence in \cite[Section III.D]{Del22a}. These two simple objects correspond to two equivalent $\CG$-module structures on $\vect \oplus \vect$ as given above.}
\end{expl}

\begin{expl}
Let $\CG = \CG(\Zb_2,\Zb_2,\alpha)$ be the finite 2-group defined as follows (recall Example \ref{expl:skeletal_2-group}):
\bit
\item Its first and second homotopy group are both $\Zb_2$.
\item The action of $\Zb_2$ on $\Zb_2$ is trivial.
\item The 3-cocycle $\alpha$ representing the associator is nontrivial. Indeed, in this case we have $H^3(\Zb_2;\Zb_2) \simeq \Zb_2$.
\eit
Let us denote $\pi_1(\CG) = \{e,x\}$ with $x^2 = e$ and $\pi_2(\CG) = \{1,y\}$ with $y^2 = 1$. The elements in $\widehat{\pi_2(\CG)}$ are denoted by $\rho_\pm$, where $\rho_\pm(y) = \pm 1$. The 3-cocycle $\alpha$ is normalized and $\alpha(x,x,x) = y$.

There are 4 simple objects in $\vect_\CG$. The simple objects $(e,\rho_+)$ and $(x,\rho_+)$ generate a fusion subcategory equivalent to $\vect_{\Zb_2}$, and the simple objects $(e,\rho_-),(x,\rho_-)$ generate a fusion subcategory equivalent to $\vect_{\Zb_2}^{(\rho_- \circ \alpha)}$. By abuse of notation, we also denote $\rho_- \circ \alpha$ by $\alpha$. Its cohomology class is also nontrivial in $H^3(\Zb_2;\bk^\times) \simeq \Zb_2$. Therefore, we have $\vect_\CG \simeq \vect_{\Zb_2} \oplus \vect_{\Zb_2}^\alpha$ as multi-fusion categories and $2\rep(\CG) \simeq 2\rep(\Zb_2) \oplus 2\rep(\Zb_2,\alpha)$ as finite semisimple 2-categories. In particular, there are 3 simple 2-representations of $\CG$.
\bit
\item The first two are two indecomposable finite semisimple $\vect_{\Zb_2}$-modules as discussed in Example \ref{expl_Z2_Z3_2-rep}. We denote them by $\one$ and $\one_c$.
\item The last one is the unique indecomposable finite semisimple $\vect_{\Zb_2}^\alpha$-module: the regular module $\vect_{\Zb_2}^\alpha$. The object $x \in \CG$ acts on $\vect_{\Zb_2}^\alpha$ by permuting two simple objects, and the morphism $y \in \pi_2(\CG)$ acts on $\vect_{\Zb_2}^\alpha$ as $\rho_-(y) = -1$. We denote this 2-representation by $T$.
\eit
This 2-category can be depicted by the following graph:
\[
\xymatrix{
\one \ar@(ul,ur)[]^{\rep(\Zb_2)}  \ar@/^/[rr]^{\vect} & & \one_c \ar@(ul,ur)[]^{\vect_{\Zb_2}} \ar@/^/[ll]^{\vect} & T \ar@(ul,ur)[]^{\vect_{\Zb_2}^\alpha}
} .
\]
Similar to the computation in Example \ref{expl_Z2_Z3_2-rep}, we obtain the fusion rules of $2\rep(\CG)$:
\[
\one_c \boxtimes \one_c \simeq \one_c \oplus \one_c , \quad T \boxtimes T \simeq \one_c \oplus \one_c , \quad \one_c \boxtimes T \simeq T \boxtimes \one_c \simeq T \oplus T .
\]
\end{expl}

\subsection{Regular algebra}

Let $\CG$ be a finite 2-group. The category $\fun(\CG,\vect)$ admits a natural symmetric monoidal structure induced by that of $\vect$. The tensor product of functors in $\fun(\CG,\vect)$ is the `point-wise' tensor product. Since $\CG$ is equivalent to the disjoint union of $\lvert \pi_1(\CG) \rvert$ copies of $\mathrm B \pi_2(\CG)$, the functor category $\fun(\CG,\vect)$ is equivalent to the direct sum of $\lvert \pi_1(\CG) \rvert$ copies of $\rep(\pi_2(\CG))$ as symmetric multi-fusion categories. Moreover, this tensor product on $\fun(\CG,\vect)$ is clearly $\CG$-equivariant if we equip $\fun(\CG,\vect)$ with the left $\CG$-module structures induced by left or right translation (see Example \ref{expl:regular_2-rep}). Therefore, $\mathcal F_l(\CG)$ and $\mathcal F_r(\CG)$ are symmetric algebras \cite{DS97,MC00} (equivalently, $E_3$-algebras or $E_\infty$-algebras) in $2\rep(\CG)$. Moreover, they are equivalent via the inverse functor $(-)^\ast \colon \CG \to \CG^\rev$. We call them the \emph{left} and \emph{right regular algebra} of $\CG$, respectively.

Suppose $\CG = \CG(G,A,\alpha)$ is skeletal. For every $g \in G$ and $\rho \in \hat A$, recall that $(g,\rho) \in \vect_\CG$ is a simple object. Using the inclusion functor $\CG \to \vect_\CG$, we can view the hom functor $\vect_\CG((g,\rho),-)$ as a functor from $\CG$ to $\vect$, denoted by $H(g,\rho) \in \fun(\CG,\vect)$. Clearly $H(g,\rho)(h) = 0$ if $g \neq h \in G$, and $H(g,\rho)(g) \simeq \bk$ is equipped with a left $\CG(g,g) = A$-action isomorphic to $\rho$. The equivalence $\vect_\CG \simeq \fun(\CG,\vect)$ in Remark \ref{rem:vect_G_dual} sends $(g,\rho) \in \vect_\CG$ to $H(g,\rho^{-1}) \in \fun(\CG,\vect)$. Thus
\[
\{H(g,\rho) \mid g \in G , \, \rho \in \hat A\}
\]
are all the simple objects of $\fun(\CG,\vect)$. We also denote $H(g) \coloneqq \vect_\CG(g,-) \in \fun(\CG,\vect)$. It sends $g \in G$ to the group algebra $\bk[A]$ and other objects to $0$.

The point-wise tensor product on $\fun(\CG,\vect)$ is given by $H(g,\rho) \otimes H(h,\sigma) \simeq \delta_{g,h} H(g,\rho\sigma)$. The left and right regular 2-representations are given by $g \odot^l H(h,\sigma) \simeq H(gh,g \triangleright \sigma)$ and $g \odot^r H(h,\sigma) \simeq H(hg^{-1},\sigma)$, respectively.

\begin{prop} \label{prop:regular_algebra_right_adjoint}
Let $\CG$ be a finite 2-group and $\omega \colon 2\rep(\CG) \to 2\vect$ be the forgetful 2-functor. Denote the right adjoint of $\omega$ by $\omega^R$. Then $\mathcal F_r(\CG) \simeq \omega^R(\vect)$ as 2-representations.
\end{prop}

\pf
For any $\CV \in 2\rep(\CG)$, a $\CG$-equivariant $\bk$-linear functor $\CV \to \mathcal F_r(\CG)$ is equivalent to a $\CG$-balanced functor $\CG \times \CV \to \vect$ that is $\bk$-linear in the second variable, where the right $\CG$-action on $\CG$ is given by the right translation. Such a $\CG$-balanced functor is clearly equivalent to a $\bk$-linear functor $\CV \to \vect$. In other words, there is an equivalence
\[
2\rep(\CG)(\CV,\mathcal F_r(\CG)) \simeq 2\vect(\omega(\CV),\vect) .
\]
Hence $\mathcal F_r(\CG) \simeq \omega^R(\vect)$.
\epf

Moreover, $2\vect$ can be viewed as a left $2\rep(\CG)$-module via $\omega$. Then the right adjoint $\omega^R$ is equivalent to the internal hom $[\vect,-] \colon 2\vect \to 2\rep(\CG)$. Hence $\omega^R(\vect) \simeq [\vect,\vect]$ is an algebra \cite[Corollary 4.2.4]{Dec24} (indeed, it is a rigid algebra by \cite[Theorem 5.2.7]{Dec24}). By tracing out the adjunction, it is not hard to see that the algebra structure on $\omega^R(\vect)$ is the regular algebra structure on $\mathcal F_r(\CG)$. Then by \cite[Theorem 5.3.4]{Dec24}, the 2-category $\RMod_{\mathcal F_r(\CG)}(2\rep(\CG))$ of right $\mathcal F_r(\CG)$-modules in $2\rep(\CG)$ is equivalent to $2\vect$ as left $2\rep(\CG)$-modules. On the other hand, by \cite[Corollary 3.12 and Lemma 3.13]{DY23}, $\RMod_{\mathcal F_r(\CG)}(2\rep(\CG))$ is a symmetric monoidal 2-category and the free module 2-functor
\[
- \otimes \mathcal F_r(\CG) \colon 2\rep(\CG) \to \RMod_{\mathcal F_r(\CG)}(2\rep(\CG))
\]
is a symmetric monoidal 2-functor. By the above discussion, the free module 2-functor can be identified with the forgetful 2-functor $\omega \colon 2\rep(\CG) \to 2\vect$.

\medskip
In the rest of this subsection, we compute the auto-equivalence 2-groups $\Aut^{\EE_3}_\CG(\mathcal F_l(\CG))$ and $\Aut^{\EE_3}_\CG(\mathcal F_r(\CG))$. Here the `auto-equivalences' both preserve the symmetric algebra structure and the $\CG$-module structure.

\begin{lem} \label{lem:symmetric_2-cocycle}
Let $A$ be a finite abelian group. Then any symmetric 2-cocycle in $Z^2(A;\bk^\times)$ is a coboundary.
\end{lem}

\pf
Given a symmetric 2-cocycle in $Z^2(A;\bk^\times)$, it suffices to show that the corresponding central extension of $A$ by $\bk^\times$ is trivial. The symmetric property implies that the extension group is also abelian, and we know that the abelian extension of $A$ by $\bk^\times$ are classified by the group $\Ext^1_\Zb(A,\bk^\times)$. Since $A$ is finite abelian, it is isomorphic to the direct sum of cyclic groups:
\[
A \simeq \bigoplus_{i=1}^m \Zb/n_i \Zb , \quad m \in \Nb , \, n_i \in \Nb_{> 0} .
\]
Therefore,
\[
\Ext^1_\Zb(A,\bk^\times) \simeq \bigoplus_{i=1}^m \Ext^1_\Zb(\Zb/n_i \Zb,\bk^\times) .
\]
But for every $n \in \Nb$ we have $\Ext^1_\Zb(\Zb/n\Zb,\bk^\times) \simeq \bk^\times / (\bk^\times)^n = 0$ because $\bk$ is algebraically closed. This completes the proof.
\epf

\begin{rem}
Indeed, a symmetric 2-cocycle in $Z^2(A;\bk^\times)$ is nothing but a 3-cocycle in the $\mathrm{E}_2$ Eilenberg-MacLane cohomology (abelian group cohomology) of $A$ \cite{EM53,EM54}. So Lemma \ref{lem:symmetric_2-cocycle} simply says that $H^3_{\mathrm{E}_2}(A;\bk^\times)$ vanishes for every finite abelian group $A$.
\end{rem}

\begin{prop} \label{prop:auto-equivalence_regular_algebra}
Let $\CG$ be a finite 2-group. The auto-equivalence 2-group $\Aut^{\EE_3}_\CG(\mathcal F_l(\CG))$ is equivalent to $\CG$, which acts on $\mathcal F_l(\CG)$ by right translation. Similarly, $\Aut^{\EE_3}_\CG(\mathcal F_r(\CG)) \simeq \CG$ acts on $\mathcal F_r(\CG)$ by left translation.
\end{prop}

\pf
Let us prove that $\Aut^{\EE_3}(\mathcal F_l(\CG)) \simeq \CG$. The proof for $\mathcal F_r(\CG)$ is similar.

For every $g \in \CG$, the right translation $g \odot^r - \colon \mathcal F_l(\CG) \to \mathcal F_l(\CG)$ is a $\CG$-equivariant symmetric monoidal equivalence. For every morphism $a \colon g \to h$, the natural transformation $a \odot^r - \colon g \odot^r \Rightarrow h \odot^r -$ is a $\CG$-equivariant symmetric monoidal natural isomorphism. Therefore, the assignment $g \mapsto g \odot^r -$ defines a monoidal functor $\CG \to \Aut^{\EE_3}_\CG(\mathcal F_l(\CG))$.

Now we assume that $\CG = \CG(G,A,\alpha)$ is skeletal. By Example \ref{expl:endomorphism_regular_2-rep}, a $\CG$-equivariant natural transformation from $(g \odot^r -)$ to $(g \odot^r -)$ is of the form $(f \odot^r -)$ for some $f \in \vect_\CG(g,g) \simeq \bk[A]$. If $(f \odot^r -)$ is a monoidal natural transformation, we have an equation of natural transformations for every $\rho,\sigma \in \hat A$:
\[
(f \odot^r -)_{H(g,\rho) \otimes H(g,\sigma)} = (f \odot^r -)_{H(g,\rho)} \otimes (f \odot^r -)_{H(g,\sigma)} .
\]
It follows that $\rho\sigma(f) = \rho(f) \sigma(f)$. So $f \in A$. This shows that the functor $\CG \to \Aut^{\EE_3}_\CG(\mathcal F_l(\CG))$ is fully faithful.

To show the essential surjectivity, we need to show that every $F \in \Aut^{\EE_3}_\CG(\mathcal F_l(\CG))$ is isomorphic to a right translation. The $\CG$-equivariance implies that there are isomorphisms
\[
F(H(gh,g \triangleright \sigma)) = F(g \odot^l H(h,\sigma)) \simeq g \odot^l F(H(h,\sigma))
\]
natural in $g \in \CG$. The naturality induces natural isomorphisms
\[
F((g,\rho) \odot^l H(h,\sigma)) \simeq (g,\rho) \odot^l H(h,\sigma) .
\]
The tensor unit of the multi-fusion category $\mathcal F_l(\CG)$ is the direct sum $\bigoplus_{g \in G} H(g,1)$ where $1 \in \hat A$ is the trivial character. As a monoidal functor, $F$ must permute these direct summands $\{H(g,1)\}_{g \in G}$. Let $g_0 \in G$ be the element satisfying $F(H(e,1)) \simeq H(g_0^{-1},1)$. Then the $\CG$-equivariance implies that
\be \label{eq_auto_regular_algebra_equivariance}
F(H(g,\rho)) = F((g,\rho) \odot^l H(e,1)) \simeq (g,\rho) \odot^l F(H(e,1)) \simeq H(gg_0^{-1},\rho) = g_0 \odot^r H(g,\rho) , \quad \forall g \in G , \, \rho \in \hat A .
\ee
So the underlying $\CG$-module functor of $F$ is naturally isomorphic to the right translation $g_0 \odot^r -$. It remains to show that the monoidal structure of $F$ coincides with that of the right translation. Note that the equivalence $F$ is the direct sum of several symmetric monoidal equivalences $F_g \colon \rep(A)_g \to \rep(A)_{gg_0^{-1}}$, where $\rep(A)_g$ is the fusion subcategory of $\mathcal F_l(\CG)$ generated by $\{H(g,\rho)\}_{\rho \in \hat A}$. By \eqref{eq_auto_regular_algebra_equivariance}, the underlying functor of $F_g$ is the identity functor on $\rep(A) \simeq \vect_{\hat A}$. A braided monoidal structure on the identity functor on $\vect_{\hat A}$ is a symmetric 2-cocycle on $\hat A$, which must be trivial by Lemma \ref{lem:symmetric_2-cocycle}. This completes the proof.
\epf

\subsection{2-group graded categories and 2-group 3-algebra}

Let $\CG$ be a finite 2-group.

\begin{defn}
A \emph{$\CG$-graded finite semisimple category} is a right $\vect_\CG$-comodule in $2\vect$.
\end{defn}

Since $\vect_\CG$ is a Hopf monoidal category, the 2-category of $\CG$-graded finite semisimple category is a monoidal 2-category. Note that $\vect_\CG$ is co-symmetric, thus $\CG$-graded finite semisimple categories can be equivalently defined as left $\vect_\CG$-comodules in $2\vect$.

Now we compute the dual Hopf monoidal category of $\vect_\CG$. Then $\CG$-graded finite semisimple categories can be equivalently defined as right (or left) modules over the dual of $\vect_\CG$. First, the dual of $\vect_\CG$ in $2\vect$ can be identified with $\fun(\CG,\vect)$:
\[
\vect_\CG^\op \simeq \fun_\bk(\vect_\CG,\vect) \simeq \fun(\CG,\vect) .
\]
The tensor product on $\fun(\CG,\vect)$ induced by the cotensor product of $\vect_\CG$ is the pointwise tensor product:
\[
(F \otimes G)(g) = (F \boxtimes G)(g \boxtimes g) = F(g) \otimes G(g) , \quad F,G \in \fun(\CG,\vect) , \, g \in \CG .
\]
The cotensor product on $\fun(\CG,\vect)$ induced by the tensor product of $\vect_\CG$ is given by
\[
\Delta(F)(g \boxtimes h) = F(g \otimes h) , \quad F \in \fun(\CG,\vect) , \, g,h \in \CG .
\]
Recall that $H(g) \in \fun(\CG,\vect)$ is the hom functor $\vect_\CG(g,-)$ for $g \in \CG$. We can compute that
\begin{multline*}
\biggl( \int^{x,y \in \CG} \vect_\CG(g,x \otimes y) \odot H(x) \boxtimes H(y)\biggr)(h \boxtimes k) = \int^{x,y \in \CG} \vect_\CG(g,x \otimes y) \otimes \vect_\CG(x,h) \otimes \vect_\CG(y,k) \\
\simeq \vect_\CG(g,h \otimes k) = H(g)(h \otimes k) , \quad g,h,k \in \CG .
\end{multline*}
Therefore, the cotensor product of $\vect_\CG$ is given by
\[
\Delta(H(g)) = \int^{x,y \in \CG} \vect_\CG(g,x \otimes y) \odot H(x) \boxtimes H(y) \simeq \int^{x \in \CG} H(x) \boxtimes H(x^L \otimes g) , \quad g \in \CG .
\]
The antipode functor on $\fun(\CG,\vect)$ is given by
\[
S(F)(g) = F(g^*) , \quad F \in \fun(\CG,\vect) , \, g \in \CG .
\]
Clearly $S(H(g)) \simeq H(g^*)$.

Let $\CV$ be a $\CG$-graded finite semisimple category. For every $g \in \CG$, the coaction of $\vect_\CG$ on $\CV$ induces a functor
\[
\CV \to \CV \boxtimes \vect_\CG \xrightarrow{1 \boxtimes H(g)} \CV \boxtimes \vect \simeq \CV .
\]
The image of $v \in \CV$ under this functor is denoted by $v_g$. Then the right $\fun(\CG,\vect)$-action on $\CV$ induced by the $\vect_\CG$-coaction is given by
\begin{align*}
\CV \boxtimes \fun(\CG,\vect) & \to \CV \\
v \boxtimes H(g) & \mapsto v_g .
\end{align*}
Note that the coevaluation 1-morphism $\vect \to \fun(\CG,\vect) \boxtimes \vect_\CG$ is given by $\bk \mapsto \int^{g \in \CG} H(g) \boxtimes g$. Thus we can also rewrite the $\vect_\CG$-coaction as follows:
\begin{align*}
\CV & \to \CV \boxtimes \vect_\CG \\
v & \mapsto \int^{g \in \CG} v_g \boxtimes g .
\end{align*}
The coassociator is a natural isomorphism
\[
\int^{g,h \in \CG} (v_g)_h \boxtimes h \boxtimes g \to \int^{g \in \CG} v_g \boxtimes g \boxtimes g , \quad v \in \CV .
\]
It induces and can be induced by a natural isomorphism
\[
(v_g)_h \simeq \vect_\CG(h,g) \odot v_g , \quad v \in \CV , \, g,h \in \CG .
\]
This can be written in a more functorial way: for every linear category $\CC$ and bilinear functor $Q \colon \CV \times \vect_\CG \to \CC$, there is a canonical isomorphism
\[
\int^{g \in \CG} Q((v_g)_h,h) \simeq Q(v_h,h) .
\]

Since $\fun(\CG,\vect)$ is a multi-fusion category, by Theorem \ref{thm:semisimple_2-cat_multi-fusion}, the 2-category $\RMod_{\fun(\CG,\vect)}(2\vect)$ of $\CG$-graded finite semisimple categories is a fusion 2-category.

\begin{expl} \label{expl_vec_G_grading}
There is an obvious right $\vect_\CG$-coaction on $\vect_\CG$ itself given by the cotensor product:
\begin{align*}
\vect_\CG & \to \vect_\CG \boxtimes \vect_\CG \\
g & \mapsto g \boxtimes g .
\end{align*}
Equivalently, the right $\fun(\CG,\vect)$-action on $\vect_\CG$ is defined by
\begin{align*}
\vect_\CG \boxtimes \fun(\CG,\vect) & \to \vect_\CG \\
g \boxtimes F & \mapsto F(g) \odot g .
\end{align*}
Thus $\vect_\CG \in \RMod_{\fun(\CG,\vect)}(2\vect)$.
\end{expl}

\begin{expl} \label{expl_graded_category_Vec}
For any $g \in \CG$, there is a right $\vect_\CG$-coaction on $\vect$ defined by
\begin{align*}
\vect & \to \vect \boxtimes \vect_\CG \\
\bk & \mapsto \bk \boxtimes g .
\end{align*}
This $\CG$-graded finite semisimple category is denoted by $\vect_g$. Equivalently, the right $\fun(\CG,\vect)$-action on $\vect_g$ is defined by
\begin{align*}
\vect \boxtimes \fun(\CG,\vect) & \to \vect \\
\bk \boxtimes F & \mapsto F(g) .
\end{align*}
When $\CG = \CG(G,A,\alpha)$ is skeletal, this module structure is induced by the monoidal functor
\[
\fun(\CG,\vect) \simeq \bigoplus_G \rep(A) \to \rep(A)_g \to \vect ,
\]
where the second functor is the projection to the $g$-th component, and the last functor is the forgetful functor. This also implies that $\End_{\fun(\CG,\vect)}(\vect_g) \simeq \End_{\rep(A)}(\vect) \simeq \vect_A$ and $\fun_{\fun(\CG,\vect)}(\vect_g,\vect_h) = 0$ for $g \neq h$.
\end{expl}

In the following, we identify $\RMod_{\fun(\CG,\vect)}(2\vect)$ with the 2-group 3-algebra of $\CG$. The \emph{2-group 3-algebra of $\CG$}, denoted by $2\vect_\CG$, is defined as follows.
\bnu
\item First, view $\CG$ as a monoidal 2-category with only identity 2-morphism and linearize this 2-category. This $\bk$-linear 2-category is denoted by $\CG \times \mathrm B^2 \bk$.
\item Take the Karoubi completion of each hom category and then take the Karoubi completion of the 2-category. So $2\vect_\CG \coloneqq \Kar(\CG \times \mathrm B^2 \bk)$.
\enu
The monoidal structure of $2\vect_\CG$ is inherited from that of $\CG$.

When $\CG = \CG(G,A,\alpha)$ is skeletal, let us compute $2\vect_\CG$ more explicitly. First, the objects of the 2-category $\CG \times \mathrm B^2 \bk$ are elements of $G$; the 1-endomorphisms on each $g \in G$ are elements of $A$ and there is no 1-morphism between different objects; the 2-endomorphism space on each 2-morphism is $\bk$ and there is no 2-morphism between different 1-morphisms. By taking the Karoubi completion of each hom category, we obtain the 2-category
\[
\coprod_G \mathrm B \vect_A .
\]
Then its Karoubi completion is
\[
\bigoplus_G \Kar(\mathrm B \vect_A) \simeq \bigoplus_G \RMod_{\vect_A}(2\vect) .
\]
This shows that $2\vect_\CG$ is finite semisimple and hence a fusion 2-category.

\begin{prop}
There is a canonical equivalence $\RMod_{\fun(\CG,\vect)}(2\vect) \simeq 2\vect_\CG$ of fusion 2-categories.
\end{prop}

\pf
It suffices to construct a monoidal 2-functor $\CG \times \mathrm B^2 \bk \to \RMod_{\fun(\CG,\vect)}(2\vect)$ and show that it is a Karoubi completion. This monoidal 2-functor is induced by a monoidal 2-functor $\CG \to \RMod_{\fun(\CG,\vect)}(2\vect)$ defined as follows.
\bit
\item An object $g \in \CG$ is mapped to $\vect_g \in \RMod_{\fun(\CG,\vect)}(2\vect)$ defined in Example \ref{expl_graded_category_Vec}.
\item A morphism $a \colon g \to g$ is mapped to the simple 1-morphism $a \in \vect_A \simeq \End_{\fun(\CG,\vect)}(\vect_g)$.
\item The monoidal structure is trivial. We need to verify the following facts about $\RMod_{\fun(\CG,\vect)}(2\vect)$.
\bnu
\item The right $\vect_\CG$-comodule structure on $\vect_g \boxtimes \vect_h$ is induced by the tensor product of $\vect_\CG$:
\[
\vect_g \boxtimes \vect_h \to \vect_g \boxtimes \vect_\CG \boxtimes \vect_h \boxtimes \vect_\CG \simeq \vect_g \boxtimes \vect_h \boxtimes \vect_\CG \boxtimes \vect_\CG \xrightarrow{1 \boxtimes 1 \boxtimes \otimes} \vect_g \boxtimes \vect_h \boxtimes \vect_\CG .
\]
It maps $\bk \boxtimes \bk$ to $\bk \boxtimes \bk \boxtimes (g \otimes h)$. Therefore, $\vect_g \boxtimes \vect_h \simeq \vect_{g \otimes h}$ as right $\vect_\CG$-comodules or right $\fun(\CG,\vect)$-modules.
\item Th associator $(\vect_g \boxtimes \vect_h) \boxtimes \vect_k \to \vect_g \boxtimes (\vect_h \boxtimes \vect_k)$ is the canonical equivalence with the right $\vect_\CG$-comodule structure induced by the associator of $\vect_\CG$. After identifying $(\vect_g \boxtimes \vect_h) \boxtimes \vect_k$ with $\vect_{(g \otimes h) \otimes k}$ and identifying $\vect_g \boxtimes (\vect_h \boxtimes \vect_k)$ with $\vect_{g \otimes (h \otimes k)}$, the associator is the identity functor equipped with the right $\vect_\CG$-comodule structure
\[
\bk \boxtimes ((g \otimes h) \otimes k) \xrightarrow{1 \boxtimes \alpha(g,h,k)} \bk \boxtimes (g \otimes (h \otimes k)) .
\]
\item The pentagonator is trivial.
\enu
\eit
To see that this monoidal 2-functor $\CG \to \RMod_{\fun(\CG,\vect)}(2\vect)$ is a Karoubi completion, it is enough to notice that
\bnu
\item $\pi_0(2\vect_\CG) \simeq G$ and $\vect_g$ is contained in the connected component corresponding to $g$ for each $g \in G$;
\item the functor $\End_{\CG \times \mathrm B^2 \bk}(g) \to \End_{\RMod_{\fun(\CG,\vect)}(2\vect)}(\vect_g) \simeq \vect_A$ is a Karoubi completion for each $g \in G$.
\enu
Hence, we conclude that the fusion 2-category $\RMod_{\fun(\CG,\vect)}(2\vect)$ of $\CG$-graded finite semisimple categories is equivalent to $2\vect_\CG$.
\epf

\begin{expl}
Suppose $\CG = G$ is a finite 1-group. Then $2\vect_G$ is the usual 2-category of $G$-graded finite semisimple categories. The tensor product of simple objects are given by the group multiplication of $G$.
\end{expl}

\begin{expl}
Suppose $\CG = \mathrm B A$ where $A$ is a finite abelian group. Then $2\vect_{\mathrm B A}$ is equivalent to $\RMod_{\vect_A}(2\vect)$ and the tensor product is given by the relative tensor product $\boxtimes_{\vect_A}$. Since $\vect_A \simeq \rep(\hat A)$ as symmetric fusion categories, $2\vect_{\mathrm B A}$ is also equivalent to $\RMod_{\rep(\hat A)}(2\vect)$ equipped with the tensor product $\boxtimes_{\rep(\hat A)}$. Then by taking the de-equivariantization \cite[Section 4]{DGNO10}, this is equivalent to $\RMod_{\vect_{\hat A}}(2\vect)$ equipped with the tensor product $\boxtimes$. In other words, $2\vect_{\mathrm B A}$ is equivalent to $2\rep(\hat A)$ as fusion 2-categories.
\end{expl}

\begin{expl}
Consider the $\CG$-graded category $\vect_\CG \in \RMod_{\fun(\CG,\vect)}(2\vect) \simeq 2\vect_\CG$ defined in Example \ref{expl_vec_G_grading}. Since $\vect_\CG$ is a Hopf monoidal category, the tensor product of $\vect_\CG$ is a $\vect_\CG$-comodule functor. Therefore, $\vect_\CG$ is also an algebra in $2\vect_\CG$. There is also an equivalent description of this algebra structure. First, for every $\CV \in 2\vect_\CG$, there is an equivalence
\begin{align*}
2\vect_\CG(2\vect)(\CV,\vect_\CG) & \simeq 2\vect(\omega{\CV},\vect) \\
F & \mapsto \bigl( \CV \xrightarrow{F} \vect_\CG \xrightarrow{\varepsilon} \vect \bigr) \\
\bigl( \CV \xrightarrow{\Delta_\CV} \CV \boxtimes \vect_\CG \xrightarrow{G \boxtimes 1} \vect \boxtimes \vect_\CG \simeq \vect_\CG \bigr) & \mapsfrom G
\end{align*}
where $\omega \colon 2\vect_\CG \to 2\vect$ is the forgetful 2-functor, $\varepsilon$ is the counit of $\vect_\CG$ and $\Delta_\CV$ is the $\vect_\CG$-coaction on $\CV$. Therefore, $\vect_\CG \simeq \omega^R(\vect)$ where $\omega^R \colon 2\vect \to 2\vect_\CG$ is the right adjoint of $\omega$. If we view $2\vect$ as a $2\vect_\CG$-module via $\omega$, the right adjoint $\omega^R$ is equivalent to the internal hom $[\vect,-] \colon 2\vect \to 2\vect_\CG$. Then the algebra structure on $\vect_\CG$ is also the same as that of the internal hom $[\vect,\vect]$.
\end{expl}

\begin{rem}
There is a notion of a ``(monoidal) category graded by a crossed module'' proposed in \cite{SV23}: given a crossed module $\chi \colon E \to H$, a $\chi$-graded category is a $\vect_E$-enriched monoidal category $\CC$ equipped with a subclass of objects $\CC_{\text{hom}}$ and a map $\CC_{\text{hom}} \to H$ satisfying certain conditions. This notion is different from the 2-group graded 2-category defined in this subsection, as explained below.
\bit
\item Although the notion of a crossed module is equivalent to a 2-group, it still encodes more information than a 2-group. For example, for any group $H$, the identity map $H \to H$ defines a crossed module $\chi$ which is equivalent to the trivial 2-group. A finite semisimple category graded by the trivial 2-group is simply a finite semisimple category, but a $\chi$-graded category in the sense of \cite{SV23} has more structures.
\item Consider the crossed module $\chi \colon A \to \ast$ where $A$ is a finite abelian group. It is equivalent to the 2-group $\mathrm B A$. A $\mathrm B A$-graded category is a module over $\fun(\mathrm B A,\vect) \simeq \rep(A) \simeq \vect_{\hat A}$, but a $\chi$-graded category is a $\vect_A$-enriched (monoidal) category (see \cite[4.3.2]{SV23}).
\eit
We expect that there is a suitably-defined notion of a $\chi$-graded finite semisimple category (without monoidal structure), and the 2-category of $\chi$-graded finite semisimple categories is equivalent to $2\vect_\CG$ where $\CG$ is a 2-group equivalent to $\chi$.
\end{rem}

\section{Tannaka-Krein duality for finite 2-groups} \label{sec:duality}

In this section, we establish the Tannaka-Krein duality for finite 2-group $\CG$. For the definition of a symmetric monoidal 2-category, a symmetric monoidal 2-functor, a symmetric monoidal 2-natural transformation and a symmetric monoidal modification, we refer readers to \cite{SP11,JY21}.

\begin{thm} \label{thm:main}
Let $\CG$ be a finite 2-group and $\omega \colon 2\rep(\CG) \to 2\vect$ be the forgetful 2-functor. Then the 2-group $\Aut^{\EE_4}(\omega)$ of symmetric monoidal 2-natural auto-equivalences of $\omega$ and symmetric monoidal modifications is equivalent to $\CG$ as 2-groups.
\end{thm}

Our proof is very similar to that for finite groups in Section \ref{sec:TK_dual_finite_group}. We prove Theorem \ref{thm:main} in three steps.

\medskip
\noindent \textbf{Step 1}: Construct a monoidal functor $\Phi \colon \CG \to \Aut^{\EE_4}(\omega)$.
\smallskip

First we construct a monoidal functor $\Phi \colon \CG \to \Aut^{\EE_4}(\omega)$. For every $g \in \CG$, we define a symmetric monoidal 2-natural auto-equivalence $\Phi^g \colon \omega \Rightarrow \omega$ as follows.
\bit
\item For every $\CV \in 2\rep(\CG)$, i.e., a finite semisimple 2-representation of $\CG$, the 1-morphism $\Phi^g_\CV \colon \omega(\CV) \to \omega(\CV)$ is defined to be the left action of $g$ on $\CV$:
\[
\Phi^g_\CV \coloneqq (g \odot -) \colon \CV \to \CV .
\]
\item For every 1-morphism $F \colon \CV \to \CW$ in $2\rep(\CG)$, i.e., a left $\CG$-module functor $F \colon \CV \to \CW$, the 2-isomorphism $\Phi^g_F$ as depicted in the following diagram
\[
\xymatrix{
\omega(\CV) \ar[r]^{\omega(F)} \ar[d]_{\Phi^g_\CV} & \omega(\CW) \ar[d]^{\Phi^g_\CW} \\
\omega(\CV) \ar[r]_{\omega(F)} & \omega(\CW) \ultwocell<\omit>{\Phi^g_F \;}
}
\]
i.e., a natural isomorphism $\Phi^g_F \colon F(g \odot -) \Rightarrow g \odot F(-)$, is defined by the left $\CG$-module structure of $F$.
\item For every $\CV,\CW \in 2\rep(\CG)$, the 2-isomorphism $\Phi^{g,2}_{\CV,\CW}$ as depicted in the following diagram
\[
\xymatrix{
\omega(\CV) \boxtimes \omega(\CW) \ar@{=}[r] \ar[d]_{\Phi^g_\CV \boxtimes \Phi^g_\CW} & \omega(\CV \boxtimes \CW) \ar[d]^{\Phi^g_{\CV \boxtimes \CW}} \\
\omega(\CV) \boxtimes \omega(\CW) \ar@{=}[r] & \omega(\CV \boxtimes \CW) \ultwocell<\omit>{\Phi^{g,2}_{\CV,\CW} \quad\;}
}
\]
i.e., a natural isomorphism $\Phi^g_\CV \boxtimes \Phi^g_\CW \Rightarrow \Phi^g_{\CV \boxtimes \CW}$, is defined to be identity.
\eit
It is not hard to see that $\Phi^g \colon \omega \Rightarrow \omega$ defined as above is a braided monoidal 2-natural equivalence (and hence a symmetric monoidal 2-natural equivalence). For every morphism $a \colon g \to h$ in $\CG$, we define a modification $\Phi^a \colon \Phi^g \Rrightarrow \Phi^h$ by
\[
\Phi^a_\CV \coloneqq \xymatrix{
\CV \rrtwocell^{\Phi^g_\CV = g \odot -}_{\Phi^h_\CV = h \odot -}{\quad a \odot -} & & \CV
} , \quad \forall \CV \in 2\rep(\CG) .
\]
It is also a monoidal modification. Then the assignment $g \mapsto \Phi^g$ and $a \mapsto \Phi^a$ defines a monoidal functor $\Phi \colon \CG \to \Aut^{\EE_4}(\omega)$.

\medskip
\noindent \textbf{Step 2}: Construct a monoidal functor $\Psi \colon \Aut^{\EE_4}(\omega) \to \CG$.
\smallskip

\begin{lem} \label{lem:symmetric_2-transformation_algebra_homomorphism}
Let $\SC,\SD$ be symmetric monoidal 2-categories, $F,G \colon \SC \to \SD$ be symmetric monoidal 2-functors, $\phi,\psi \colon F \Rightarrow G$ be symmetric monoidal 2-natural transformations and $\Gamma \colon \phi \Rrightarrow \psi$ be a symmetric monoidal modification. Then for any symmetric algebra $A \in \SC$, the component $\phi_A \colon F(A) \to G(A)$ defines a symmetric algebra 1-homomorphism and $\Gamma_A \colon \phi_A \Rightarrow \psi_A$ is a symmetric algebra 2-homomorphism.
\end{lem}

\pf
We view the symmetric algebra $A \in \SC$ as a symmetric lax monoidal 2-functor (also called a symmetric weak monoidal homomorphism in \cite{DS97}) $A \colon \ast \to \SC$, where $\ast$ is the trivial 2-category with only one object, identity 1-morphism and identity 2-morphism. Then the component $\phi_A$ can be viewed as the whiskering $\phi \circ A$, which is a symmetric monoidal 2-natural transformation from $F \circ A = F(A)$ to $G \circ A = G(A)$. Similarly, $\Gamma_A$ is also the whiskering $\Gamma \circ A$.
\epf

\begin{rem}
It is also useful to write down all the defining data of the symmetric algebra 1-homomorphism $\phi_A \colon F(A) \to G(A)$.
\bnu[(1)]
\item There is a 2-isomorphism defined by the composition of the following two 2-morphisms in $\SD$:
\[
\xymatrix{
F(A) \otimes F(A) \ar[r]^-{F^2_{A,A}} \ar[d]_{\phi_A \otimes \phi_A} & F(A \otimes A) \ar[r]^-{F(\mu)} \ar[d]_{\phi_{A \otimes A}} & F(A) \ar[d]^{\phi_A} \\
G(A) \otimes G(A) \ar[r]_-{G^2_{A,A}} & G(A \otimes A) \ar[r]_-{G(\mu)} \ultwocell<\omit>{\phi^2_{A,A} \quad} & G(A) \ultwocell<\omit>{\phi_\mu \;\;}
}
\]
where $\mu \colon A \otimes A$ is the multiplication.
\item There is a 2-isomorphism defined by the composition of the following two 2-morphisms in $\SD$:
\[
\xymatrix{
\one \ar[r]^{F^0} \ar[dr]_{G^0} & F(\one) \ar[r]^{F(\eta)} \ar[d]^{\phi_\one} & F(A) \ar[d]^{\phi_A} \\
 & G(\one) \ar[r]_{G(\eta)} \ultwocell<\omit>{<2>\;\phi^0} & G(A) \ultwocell<\omit>{\phi_\eta \; \;}
}
\]
where $\eta \colon \one \to A$ is the unit of the algebra $A$.
\enu
\end{rem}

Let $\phi \in \Aut^{\EE_4}(\omega)$. For any $\CV \in 2\rep(\CG)$, the naturality of $\phi$ means that $\phi_\CV \colon \omega(\CV) \to \omega(\CV)$ is equipped with 2-morphisms $\{\phi_f\}$ for all $f \in \End_{2\rep(\CG)}(\CV)$, as depicted in the following diagram:
\[
\xymatrix{
\omega(\CV) \ar[r]^-{\phi_\CV} \ar[d]_{\omega(f)} \drtwocell<\omit>{\;\phi_f} & \omega(\CV) \ar[d]^{\omega(f)} \\
\omega(\CV) \ar[r]^-{\phi_\CV} & \omega(\CV)
}
\]
The axioms of a 2-natural transformation imply that $\{\phi_f\}$ is a half-braiding on $\phi_\CV$. In other words, $\phi_\CV$ can be viewed as an object in the Drinfeld centralizer
\[
\FZ_1 \bigl(\End_{2\rep(\CG)}(\CV) \xrightarrow{\omega} \End_{2\vect}(\omega(\CV)) \bigr) \simeq \End_{\End_{2\rep(\CG)}(\CV)}(\omega(\CV)) .
\]
In particular, $\phi_{\mathcal F_r(\CG)} \in \End_{\vect_\CG}(\omega(\mathcal F_r(\CG))) \simeq \End_{2\rep(\CG)}(\mathcal F_l(\CG))$ is $\CG$-equivariant with respect to the left translation $\CG$-action (see Example \ref{expl:endomorphism_regular_2-rep}). By Lemma \ref{lem:symmetric_2-transformation_algebra_homomorphism}, the restriction $\phi \mapsto \phi_{\mathcal F_r(\CG)}$ defines a monoidal functor $\Aut^{\EE_4}(\omega) \to \Aut^{\EE_3}_\CG(\mathcal F_l(\CG))$. Then by Proposition \ref{prop:auto-equivalence_regular_algebra} we obtain a monoidal functor
\[
\Psi \colon \Aut^{\EE_4}(\omega) \to \Aut^{\EE_3}_\CG(\mathcal F_l(\CG)) \simeq \CG .
\]

\medskip
\noindent \textbf{Step 3}: Show that $\Phi$ and $\Psi$ are quasi-inverses to each other.
\smallskip

It is easy to see that $\Psi \circ \Phi \colon \CG \to \CG$ is isomorphic to identity. It remains to show that $\Phi \circ \Psi \colon \Aut^{\EE_4}(\omega) \to \Aut^{\EE_4}(\omega)$ is isomorphic to identity. In other words, we need to show that any $\phi \in \Aut^{\EE_4}(\omega)$ with $\phi_{\mathcal F_r(\CG)} = g \odot^r -$ is isomorphic to $\Phi^g$.

Note that $\mathcal F_r(\CG) \simeq \vect_G \in 2\rep(\CG)$ is a generator, i.e., for every $\CV \in 2\rep(\CG)$ the hom category between $\CV$ and $\mathcal F_r(\CG)$ is nonzero (see Proposition \ref{prop:regular_algebra_right_adjoint}). Therefore, the natural inclusion $\iota \colon \mathrm B \End_{2\rep(\CG)}(\mathcal F_r(\CG)) \to 2\rep(\CG)$ is a Karoubi completion \cite{DR18} or a condensation completion \cite{GJF19}. By the universal property of the Karoubi completion, the restriction $\phi \mapsto \phi_{\mathcal F_r(\CG)}$ defines an equivalence (see also Section \ref{sec:2Rep(G)_2-cat})
\begin{multline*}
\End(\omega) \to \End(\omega \circ \iota) \simeq \FZ_1 \bigl( \End_{2\rep(\CG)}(\mathcal F_r(\CG)) \xrightarrow{\omega} \End_{2\vect}(\omega(\mathcal F_r(\CG))) \bigr) \\
\simeq \End_{\vect_\CG}(\fun(\CG,\vect)) \simeq \End_{2\rep(\CG)}(\mathcal F_l(\CG)) \simeq \vect_\CG .
\end{multline*}
If $\phi_{\mathcal F_r(\CG)} = g \odot^r -$, then by the above result we see that the underlying 2-natural equivalence of $\phi$ must be $\Phi^g$.

Now we only need to consider the monoidal structure of $\phi$, which is an invertible modification $\phi^2$ whose components are depicted in the following diagram:
\[
\xymatrix{
\omega(\CV) \boxtimes \omega(\CW) \ar@{=}[r] \ar[d]_{\phi_\CV \boxtimes \phi_\CW} & \omega(\CV \boxtimes \CW) \ar[d]^{\phi_{\CV \boxtimes \CW}} \\
\omega(\CV) \boxtimes \omega(\CW) \ar@{=}[r] & \omega(\CV \boxtimes \CW) \ultwocell<\omit>{\phi^2_{\CV,\CW} \quad\;}
}
\]
Using the same Karoubi completion argument, $\phi^2$ is completely determined by the component $\phi^2_{\vect_\CG,\vect_\CG}$ because $\vect_\CG \in 2\rep(\CG)$ is a generator. We can also give a more explicit construction: for every $\CV,\CW \in 2\rep(\CG)$ and $v \in \CV$, $w \in \CW$, since $- \odot v \colon \vect_\CG \to \CV$ and $- \odot w \colon \vect_\CG \to \CW$ are 1-morphisms in $2\rep(\CG)$, the modification axiom for $\phi^2$ implies that the following diagram commutes:
\[
\xymatrix@C=7em{
(g \boxtimes g) \odot (v \boxtimes w) \ar[r]^-{(\phi^2_{\vect_\CG,\vect_\CG})_{e,e} \odot 1} \ar@{=}[d] & (g \boxtimes g) \odot (v \boxtimes w) \ar@{=}[d] \\
(g \odot v) \boxtimes (g \odot w) \ar[r]^-{(\phi^2_{\CV,\CW})_{v,w}} & (g \odot v) \boxtimes (g \odot w)
}
\]
Thus $\phi^2$ is completely determined by the morphism $(\phi^2_{\vect_\CG,\vect_\CG})_{e,e} \colon g \boxtimes g \to g \boxtimes g$. Decomposing $e$ as the direct sum of simple objects $e = \bigoplus_{\rho \in \hat A} (e,\rho)$, the morphism $(\phi^2_{\vect_\CG,\vect_\CG})_{e,e}$ is the direct sum of
\[
(\phi^2_{\vect_\CG,\vect_\CG})_{(e,\rho),(e,\sigma)} \colon (g \otimes (e,\rho)) \boxtimes (g \otimes (e,\sigma)) = (g,g \triangleright \rho) \boxtimes (g,g \triangleright \sigma) \to (g,g \triangleright \rho) \boxtimes (g,g \triangleright \sigma)
\]
for all $\rho,\sigma \in \hat A$. We denote
\[
\tilde \phi_{\rho,\sigma} \coloneqq (\phi^2_{\vect_\CG,\vect_\CG})_{(e,\rho),(e,\sigma)} \in \Aut_{\vect_\CG \boxtimes \vect_\CG}((g,g \triangleright \rho) \boxtimes (g,g \triangleright \sigma)) \simeq \bk^\times .
\]
Then $\tilde \phi \in C^2(\hat A;\bk^\times)$. The monoidal axiom for $\phi^2$ implies that $\tilde \phi \in Z^2(\hat A;\bk^\times)$ is a 2-cocycle, and the braided axiom for $\phi^2$ implies that $\tilde \phi_{\rho,\sigma} = \tilde \phi_{\sigma,\rho}$ is symmetric.

By Lemma \ref{lem:symmetric_2-cocycle}, $\tilde \phi = \mathrm d \tilde \Gamma$ for some $\tilde \Gamma \in C^1(\hat A;\bk^\times)$. Then the direct sum of morphisms
\[
(\Gamma_{\vect_\CG})_e \coloneqq \biggl(g \otimes e = \bigoplus_{\rho \in \hat A} (g,g \triangleright \rho) \xrightarrow{\bigoplus \tilde \Gamma_{\rho}} \bigoplus_{\rho \in \hat A} (g,g \triangleright \rho) = g \otimes e \biggr)
\]
induces a modification $\Gamma \colon \phi \Rrightarrow \Phi^g$ because $\vect_\CG \in 2\rep(\CG)$ is a generator. More explicitly, for any $\CV \in 2\rep(\CG)$ and $v \in \CV$, the modification axiom implies that
\[
(\Gamma_\CV)_v = (\Gamma_{\vect_\CG})_e \odot 1_v \colon g \odot v \to g \odot v .
\]
The equation $\tilde \phi = \mathrm d \tilde \Gamma$ is equivalent to that $\Gamma \colon \phi \Rrightarrow \Phi^g$ is a monoidal modification. Hence we prove that $\phi$ is isomorphic to $\Phi^g$ as symmetric monoidal 2-natural equivalences, and the proof of Theorem \ref{thm:main} is complete.

\bibliography{Top}

\begin{thebibliography}{BBFW12}

\bibitem[BBFW12]{BBFW12}
John Baez, Aristide Baratin, Laurent Freidel, and Derek Wise.
\newblock {\em Infinite-dimensional representations of 2-groups}.
\newblock Memoirs of the American Mathematical Society, 219(1032):0--0, 2012.
\newblock \href{https://arxiv.org/abs/0812.4969}{arXiv:0812.4969}.

\bibitem[BC04]{BC04}
John~C. Baez and Alissa~S. Crans.
\newblock {\em Higher-dimensional algebra {VI}: {Lie} 2-algebras}.
\newblock Theory and Applications of Categories, 12(15):492--528, 2004.
\newblock \href{https://arxiv.org/abs/math/0307263}{arXiv:math/0307263}.

\bibitem[BL04]{BL04}
John~C. Baez and Aaron~D. Lauda.
\newblock {\em Higher-dimensional algebra {V}: 2-groups}.
\newblock Theory and Applications of Categories, 12(14):423--491, 2004.
\newblock \href{https://arxiv.org/abs/math/0307200}{arXiv:math/0307200}.

\bibitem[BM06]{BM06}
John~W. Barrett and Marco Mackaay.
\newblock {\em Categorical representations of categorical groups}.
\newblock Theory and Applications of Categories, 16(20):529--557, 2006.
\newblock \href{https://arxiv.org/abs/math/0407463}{arXiv:math/0407463}.

\bibitem[CF94]{CF94}
Louis Crane and Igor~B. Frenkel.
\newblock {\em Four-dimensional topological quantum field theory, {Hopf}
  categories, and the canonical bases}.
\newblock Journal of Mathematical Physics, 35(10):5136--5154, 1994.
\newblock \href{https://arxiv.org/abs/hep-th/9405183}{arXiv:hep-th/9405183}.

\bibitem[CY05]{CY05}
Louis Crane and David~N. Yetter.
\newblock {\em Measurable categories and 2-groups}.
\newblock Applied Categorical Structures, 13(5-6):501--516, 2005.
\newblock \href{https://arxiv.org/abs/math/0305176}{arXiv:math/0305176}.

\bibitem[D{\'{e}}c23]{Dec23}
Thibault D{\'{e}}coppet.
\newblock {\em Compact semisimple 2-categories}.
\newblock Transactions of the American Mathematical Society, 2023.
\newblock \href{https://arxiv.org/abs/2111.09080}{arXiv:2111.09080}.

\bibitem[D{\'{e}}c24]{Dec24}
Thibault~D. D{\'{e}}coppet.
\newblock {\em Finite semisimple module 2-categories}.
\newblock Selecta Mathematica, 31(1), 2024.
\newblock \href{https://arxiv.org/abs/2107.11037}{arXiv:2107.11037}.

\bibitem[Del02]{Del02}
P.~Deligne.
\newblock {\em Cat\'{e}gories tensorielles}.
\newblock Moscow Mathematical Journal, 2(2):227--248, 2002.

\bibitem[Del07]{Del07}
P.~Deligne.
\newblock {\em Cat{\'{e}}gories tannakiennes}.
\newblock In {\em The Grothendieck Festschrift}, pages 111--195. Birkhäuser
  Boston, 2007.

\bibitem[Del22]{Del22a}
Clement Delcamp.
\newblock {\em A toy model for categorical charges}.
\newblock arXiv preprint, 2022.
\newblock \href{https://arxiv.org/abs/2208.07361}{arXiv:2208.07361}.

\bibitem[DGNO10]{DGNO10}
Vladimir Drinfeld, Shlomo Gelaki, Dmitri Nikshych, and Victor Ostrik.
\newblock {\em On braided fusion categories {I}}.
\newblock Selecta Mathematica, 16(1):1--119, 2010.
\newblock \href{https://arxiv.org/abs/0906.0620}{arXiv:0906.0620}.

\bibitem[DR18]{DR18}
Christopher~L. Douglas and David~J. Reutter.
\newblock {\em Fusion 2-categories and a state-sum invariant for 4-manifolds}.
\newblock arXiv preprint, 2018.
\newblock \href{https://arxiv.org/abs/1812.11933}{arXiv:1812.11933}.

\bibitem[DS97]{DS97}
Brian Day and Ross Street.
\newblock {\em Monoidal bicategories and {Hopf} algebroids}.
\newblock Advances in Mathematics, 129(1):99--157, 1997.

\bibitem[DY23]{DY23}
Thibault~D. D{\'{e}}coppet and Matthew Yu.
\newblock {\em Gauging noninvertible defects: a 2-categorical perspective}.
\newblock Letters in Mathematical Physics, 113(2), 2023.
\newblock \href{https://arxiv.org/abs/2211.08436}{arXiv:2211.08436}.

\bibitem[EGNO15]{EGNO15}
Pavel Etingof, Shlomo Gelaki, Dmitri Nikshych, and Victor Ostrik.
\newblock {\em Tensor categories}.
\newblock American Math. Soc, Providence, RI, 2015.

\bibitem[Elg07]{Elg07}
Josep Elgueta.
\newblock {\em Representation theory of 2-groups on {Kapranov} and
  {Voevodsky}‘s 2-vector spaces}.
\newblock Advances in Mathematics, 213(1):53--92, 2007.
\newblock \href{https://arxiv.org/abs/math/0408120}{arXiv:math/0408120}.

\bibitem[Elg11]{Elg11}
Josep Elgueta.
\newblock {\em On the regular representation of an (essentially) finite
  2-group}.
\newblock Advances in Mathematics, 227(1):170--209, 2011.
\newblock \href{https://arxiv.org/abs/0907.0978}{arXiv:0907.0978}.

\bibitem[EM53]{EM53}
Samuel Eilenberg and Saunders MacLane.
\newblock {\em On the groups {$H(\Pi, n)$}, {I}}.
\newblock The Annals of Mathematics, 58(1):55, 1953.

\bibitem[EM54]{EM54}
Samuel Eilenberg and Saunders MacLane.
\newblock {\em On the groups {$H(\Pi, n)$}, {II}: Methods of computation}.
\newblock The Annals of Mathematics, 60(1):49, 1954.

\bibitem[ENO10]{ENO10}
Pavel Etingof, Dmitri Nikshych, and Viktor Ostrik.
\newblock {\em Fusion categories and homotopy theory}.
\newblock Quantum Topology, pages 209--273, 2010.
\newblock \href{https://arxiv.org/abs/0909.3140}{arXiv:0909.3140}.

\bibitem[FB03]{FB03}
Magnus Forrester-Barker.
\newblock {\em Representations of crossed modules and {Cat}$^1$-groups}.
\newblock Ph.D. thesis, 2003.

\bibitem[GJF19]{GJF19}
Davide Gaiotto and Theo Johnson-Freyd.
\newblock {\em Condensations in higher categories}.
\newblock arXiv preprint, 2019.
\newblock \href{https://arxiv.org/abs/1905.09566}{arXiv:1905.09566}.

\bibitem[Gre23]{Gre23}
David Green.
\newblock {\em Tannaka-{Krein} reconstruction for fusion 2-categories}.
\newblock arXiv preprint, 2023.
\newblock \href{https://arxiv.org/abs/2309.05591}{arXiv:2309.05591}.

\bibitem[HE16]{HE16}
Benjamin~A. Heredia and Josep Elgueta.
\newblock {\em On the representations of 2-groups in {Baez}-{Crans} 2-vector
  spaces}.
\newblock Theory and Applications of Categories, 31(32):907--927, 2016.
\newblock \href{https://arxiv.org/abs/1607.04986}{arXiv:1607.04986}.

\bibitem[HXZ24]{HXZ24}
Mo~Huang, Hao Xu, and Zhi-Hao Zhang.
\newblock {\em The 2-character theory of finite 2-groups}.
\newblock arXiv preprint, 2024.
\newblock \href{https://arxiv.org/abs/2404.01162}{arXiv:2404.01162}.

\bibitem[JS91]{JS91a}
Andr{\'{e}} Joyal and Ross Street.
\newblock {\em An introduction to {Tannaka} duality and quantum groups}.
\newblock In {\em Lecture Notes in Mathematics}, pages 413--492. Springer
  Berlin Heidelberg, 1991.

\bibitem[JS93]{JS93}
A.~Joyal and R.~Street.
\newblock {\em Braided tensor categories}.
\newblock Advances in Mathematics, 102(1):20--78, 1993.

\bibitem[JY21]{JY21}
Niles Johnson and Donald Yau.
\newblock {\em 2-Dimensional Categories}.
\newblock Oxford University Press, 2021.
\newblock \href{https://arxiv.org/abs/2002.06055}{arXiv:2002.06055}.

\bibitem[Kre49]{Kre49}
Mark~Grigorievich Krein.
\newblock {\em A principle of duality for bicompact groups and quadratic block
  algebras}.
\newblock Doklady Akademii Nauk SSSR, 69:725--728, 1949.

\bibitem[KV94]{KV94}
M.~M. Kapranov and V.~A. Voevodsky.
\newblock {\em 2-categories and {Zamolodchikov} tetrahedra equations}.
\newblock In {\em Algebraic Groups and Their Generalizations: Quantum and
  Infinite-Dimensional Methods}, volume~56 of {\em Proceedings of Symposia in
  Pure Mathematics}, pages 177--260. American Mathematical Society, 1994.

\bibitem[KZ18]{KZ18}
Liang Kong and Hao Zheng.
\newblock {\em The center functor is fully faithful}.
\newblock Advances in Mathematics, 339:749--779, 2018.
\newblock \href{https://arxiv.org/abs/1507.00503}{arXiv:1507.00503}.

\bibitem[KZ24]{KZ24}
Liang Kong and Hao Zheng.
\newblock {\em Categories of quantum liquids {II}}.
\newblock Communications in Mathematical Physics, 405(9), 2024.
\newblock \href{https://arxiv.org/abs/2107.03858}{arXiv:2107.03858}.

\bibitem[Mac63]{Mac63}
Saunders MacLane.
\newblock {\em Natural associativity and commutativity}.
\newblock Rice University Studies, 49(4):28--46, 1963.

\bibitem[Mas98]{Mas98}
Heinrich Maschke.
\newblock {\em Ueber den arithmetischen charakter der coefficienten der
  substitutionen endlicher linearer substitutionsgruppen}.
\newblock Mathematische Annalen, 50(4):492--498, 1898.

\bibitem[Mas99]{Mas99}
Heinrich Maschke.
\newblock {\em Beweis des satzes, dass diejenigen endlichen linearen
  substitutionsgruppen, in welchen einige durchgehends verschwindende
  coefficienten auftreten, intransitiv sind}.
\newblock Mathematische Annalen, 52(2-3):363--368, 1899.

\bibitem[McC00]{MC00}
Paddy McCrudden.
\newblock {\em Balanced coalgebroids}.
\newblock Theory and Applications of Categories, 7(6):71--147, 2000.

\bibitem[Neu97]{Neu97}
Martin Neuchl.
\newblock {\em Representation theory of {Hopf} categories}.
\newblock Doctoral dissertation, 1997.

\bibitem[S{\'{i}}n75]{Sin75}
Hoàng~Xuân S{\'{i}}nh.
\newblock {\em Gr-catégories}.
\newblock Ph.D. thesis, 1975.

\bibitem[SP11]{SP11}
Christopher~J. Schommer-Pries.
\newblock {\em The classification of two-dimensional extended topological field
  theories}.
\newblock arXiv preprint, 2011.
\newblock \href{https://arxiv.org/abs/1112.1000}{arXiv:1112.1000}.

\bibitem[SV23]{SV23}
Kürşat Sözer and Alexis Virelizier.
\newblock {\em Monoidal categories graded by crossed modules and 3-dimensional
  hqfts}.
\newblock Advances in Mathematics, 428:109155, 2023.
\newblock \href{https://arxiv.org/abs/2207.06534}{arXiv:2207.06534}.

\bibitem[Tan39]{Tan39}
Tadao Tannaka.
\newblock {\em Über den dualitätssatz der nichtkommutativen topologischen
  gruppen}.
\newblock Tohoku Mathematical Journal, First Series, 45:1--12, 1939.

\bibitem[Yet05]{Yet05}
D.~N. Yetter.
\newblock {\em Measurable categories}.
\newblock Applied Categorical Structures, 13(5-6):469--500, 2005.
\newblock \href{https://arxiv.org/abs/math/0309185}{arXiv:math/0309185}.

\end{thebibliography}

\end{document}